\documentclass[final,1p,times,singlecolumn]{elsarticle}

\usepackage{color,fontenc,float}
\usepackage{tikz}
\usepackage[normalem]{ulem}
\usepackage{amsmath,amsthm,amssymb,framed}
\journal{Physica D}
\usepackage{mathtools}
\biboptions{sort&compress}

\newcommand{\rem}[1]{}
\newcommand{\bom}{\mbox{\boldmath$\omega$}}

\newcommand{\bu}{\mbox{\boldmath$u$}}
\newcommand{\bk}{\mbox{\boldmath$k$}}

\newcommand{\bx}{\mbox{\boldmath$x$}}

\newcommand{\I}{\int_{V}}
\newcommand{\Rn}{{\rm Re}_{\nu}}
\newcommand{\Rb}{{\rm Re}_{\beta}}
\newcommand{\Ra}{\mathcal{A}_{0}}
\newcommand{\bel}{\begin{equation}\label}
\newcommand{\ee}{\end{equation}}
\newcommand{\beq}{\begin{eqnarray}\label} 
\newcommand{\eeq}{\end{eqnarray}} 
\newcommand{\bc}{\begin{center}} 
\newcommand{\ec}{\end{center}} 
\newcommand{\ben}{\begin{enumerate}}
\newcommand{\een}{\end{enumerate}}
\newcommand{\bit}{\begin{itemize}}
\newcommand{\eit}{\end{itemize}}
\newcommand\shalf{\ensuremath{{\scriptstyle\frac{1}{2}}}}
\newcommand\sthird{\ensuremath{{\scriptstyle\frac{1}{3}}}}
\newcommand\twothirds{\ensuremath{{\scriptstyle\frac{2}{3}}}}
\newcommand\squart{\ensuremath{{\scriptstyle\frac{1}{4}}}}
\newcommand\sthreequart{\ensuremath{{\scriptstyle\frac{3}{4}}}}

\theoremstyle{definition}

\numberwithin{equation}{section}

\begin{document}
\begin{frontmatter}
\title{An analytical and computational study of the incompressible Toner-Tu Equations}

\author[inst1]{John D. Gibbon}
\ead{j.d.gibbon@ic.ac.uk}
\address[inst1]{Department of Mathematics, Imperial College London, London SW7 2AZ, UK}

\author[inst2]{Kolluru Venkata Kiran}
\ead{kollurukiran@iisc.ac.in}
\address[inst2]{Centre for Condensed Matter Theory, Department of Physics, Indian Institute of Science, Bangalore, 560012, India}

\author[inst2]{Nadia Bihari Padhan}
\ead{nadia@iisc.ac.in}

\author[inst2]{Rahul Pandit}
\ead{rahul@iisc.ac.in}

\begin{abstract}

The incompressible Toner-Tu (ITT) partial differential equations (PDEs) are an important example of a set of active-fluid PDEs. While they share certain properties with the Navier-Stokes equations (NSEs), such as the same scaling invariance, there are also important differences. The NSEs are usually considered in either the decaying or the additively forced cases, whereas the ITT equations have no additive forcing. Instead, they include a linear, activity term $\alpha \bu$ ($\bu$ is the velocity field) which pumps energy into the system, but also a negative $\bu|\bu|^{2}$-term which provides a platform for either frozen or statistically steady states. Taken together, these differences make the ITT equations an intriguing candidate for study using a combination of PDE analysis and pseudo-spectral direct numerical simulations (DNSs). In the $d=2$ case, we have established global regularity of solutions, but we have also shown the existence of bounded hierarchies of weighted, time-averaged norms of both higher derivatives and higher moments of the velocity field. Similar bounded hierarchies for Leray-type weak solutions have also been established in the $d=3$ case. We present results for these norms from our DNSs in both $d=2$ and $d=3$, and contrast them with their Navier-Stokes counterparts.
\end{abstract}
\date{\today} 

\end{frontmatter}

\section{Introduction}

The field of \textit{active matter} continues to grow rapidly \cite{R1987,VCBCS1995,CSV1997,TT1995,TTR2005,AJC2020,WDHDGLY2012,DHDWBG2013,SD2015,LMBE2019,LMBE2020,BJF2015,
ACJ2021,MJRLP2013,G2015,SJ2020,SD2017,MSJR2021,KGVP2022,SH1977,RP2022,RP2020,CLT2016,CLT2018,BFMR2022}. The term is generally used for systems that have bodies, e.g. \textit{birdoids} in computer animations \cite{R1987}, birds in a flock~\cite{VCBCS1995,TT1995,TTR2005}, or bacteria in dense suspensions \cite{TT1995,TTR2005,AJC2020,WDHDGLY2012,DHDWBG2013,SD2015,LMBE2019,LMBE2020,BJF2015,ACJ2021,SD2017,MJRLP2013,G2015,SJ2020}, all of which use some source of energy, typically internal, to move or to apply forces.  Such bodies, referred to as \textit{active particles} in the physics literature, mutually interact and lead to non-equilibrium states, which may display rich spatio-temporal evolution. The bird-flocking model of Vicsek~\cite{VCBCS1995}, a non-equilibrium version of a Heisenberg-spin model, is defined in discrete time for an assembly of point particles which are distributed randomly in space. These particles try to align with their neighbours, but with some stochastically modelled  error.

Soon after the publication of the Vicsek model \cite{VCBCS1995}, Toner and Tu (TT) introduced a hydrodynamic stochastic partial differential equation (PDE) that models flocking phenomena~\cite{TT1995,TTR2005}. The TT velocity field obeys a generalised, compressible, stochastically forced set of Navier-Stokes equations (NSEs), which are not Galilean invariant. Other hydrodynamic PDEs, related to the NSEs, were then developed to study the spatio-temporal evolution of active fluids, such as dense bacterial suspensions
\cite{TT1995,TTR2005,AJC2020,WDHDGLY2012,DHDWBG2013,SD2015,LMBE2019,LMBE2020,BJF2015,ACJ2021,SD2017,MJRLP2013,G2015,SJ2020}, or \textit{active nematics}. In one of the simplest variants, called either the mean-bacterial-velocity or the Toner-Tu-Swift-Hohenberg (TTSH) model, a term consisting of the sum of a negative Laplacian and a bi-Laplacian is added to an incompressible, deterministic TT PDE (henceforth, ITT)~\cite{SH1977}. For recent studies of the stochastically forced and deterministic variants of the ITT equations we refer the reader to Refs. \cite{RP2020,RP2022,CLT2016,CLT2018,BFMR2022}.

Although these active-matter and active-fluid PDEs have been studied intensively over the past two decades from a physical perspective, with the results of these investigations having undergone wide experimental comparison, detailed methods of Navier-Stokes analysis \cite{Leray1934,FMRT2001,DG1995,RRS2016,FGT1981,JDG2019,JDG2020} have not generally been applied directly to the ITT equations\footnote{An exception is the work of Zanger, L\"owen, and Saal on the regularity of solutions of the TTSH equations \cite{ZLS2015}.}. However, models related to the NSEs with an absorption term have been studied \cite{Anton2010}, including the Brinkman-Forchheimer-extended Darcy model of porous media~
\cite{TT2022,Anton2010,Cai2008,Celebi2006,Liu2007,LSST2015,LST2015,MTT2016,Payne1999,BS2008,Ugurlu2008,Wang2008,You2012}. The major difference is that these models possess a nonlinearity that breaks the Navier-Stokes invariance enjoyed by the ITT equations. 

The main aim of this paper is to consider the ITT PDEs in $d$ spatial dimensions ($d=2,3$) using the ideas developed in Refs.~\cite{DKGGPV2013,GDKGPV2014,GGKPPPPSS2016,GGPP2018} in which a combination of analysis and direct numerical simulations (DNSs) on the $d=3$ NSEs were used to match the results of the former against those of the latter. For the NSEs in dimension $d=3$, only weak solutions (in the sense of Leray \cite{Leray1934}) are known to exist. To elevate these to the status of strong solutions, a uniqueness property would need to be proved, but this proof remains elusive \cite{Leray1934,FMRT2001,DG1995,RRS2016,FGT1981,JDG2019,JDG2020}. The first task is to show that the ITT equations in dimension $d=3$ have equivalent weak solution properties to their Navier-Stokes counterpart. The second task is to study the $d=2$ case, keeping in mind that the NSEs are known to be regular, i.e. global existence and uniqueness of solutions have been proved \cite{FMRT2001}. It will be shown that the ITT equations also turn out to be regular, although the counterpart of bounds that remain fixed for the NSEs actually grow exponentially in time. While numerical (pseudo-spectral) methods and experiments are able to track a solution that evolves from \textit{specified} initial conditions, methods of analysis are unable to do this; instead, in a complementary fashion, they provide us with constraints on solutions that evolve from \textit{all smooth initial conditions}. Thus, one should not expect the estimated bounds to be saturated, as these take into account all smooth initial conditions, however large -- see ~\cite{DKGGPV2013,GDKGPV2014,GGKPPPPSS2016,GGPP2018}. Methods of analysis also provide upper bounds on average inverse length scales, which can be interpreted as lower bounds on the grid sizes necessary to resolve solutions.

In \S\ref{sec:dimless} we define the PDEs in dimensionless form and the quantities that are required for our analysis. In \S\ref{sec:scaling} the scaling properties of the Navier-Stokes and the ITT equations are discussed and how their similarity acts as a guide to our choice of moments of higher derivatives of the velocity field. In \S\ref{sec:energy} we discuss energy estimates. In \S\ref{sec:num} we describe the pseudo-spectral DNS that has been used to solve the ITT equations. \S\ref{sec:sumd2} deals solely with the $d=2$ case\,: in \S\ref{subsec:est2} we present a summary of our results, the proofs of which can be found in \ref{appA},  while \S\ref{subsec:reg1} and \S\ref{subsec:num2} are devoted to global regularity and numerical results respectively. Likewise in \S\ref{sec:sumd3} we summarise our results for $d=3$, with the proofs found in \ref{appB}. In \S\ref{sec:Conclusions} we discuss the significance of our results and compare them with similar results for related PDEs.


\section{Dimensionless equations}\label{sec:dimless}

The standard form of the incompressible Toner-Tu (ITT) equations is given by \cite{TT1995,TTR2005,RP2020}\,: 
\bel{a1a}
\left(\partial_{t}+\lambda\bu\cdot\nabla\right)\bu + \nabla p = \alpha\bu + \nu\Delta\bu - \beta\bu|\bu|^{2}\,.
\ee
The fixed parameters $\alpha,\,\beta$ are positive and the velocity field $\bu$ satisfies the incompressibility condition $\mbox{div}\,\bu = 0$. $\beta$ has the dimension $TL^{-2} \equiv [\nu^{-1}]$, $\alpha$ is a frequency and $\lambda$ is a dimensionless parameter. The domain is taken to be a periodic box $[0,\,L]^{d}_{per}$. We leave remarks until \S\ref{sec:scaling} on the literature involving generalizations of this system to a $\bu|\bu|^{2\delta}$ nonlinear term.

The first step is to introduce a typical velocity field $U_{0}$ for which we have two definitions:
\bel{a1b}
U_{0} = \sqrt{\alpha/\beta}\qquad\mbox{and}\qquad U_{0} = \nu/L\,.
\ee
Then primed dimensionless variables are defined thus\,: 
\beq{nv1}
\bx' = L^{-1}\bx\,; \qquad t' = U_{0}L^{-1}t\,;\qquad
\bu' = \lambda U_{0}^{-1}\bu\,; \qquad p' = \lambda U_{0}^{-2}p\,.
\eeq
This transforms (\ref{a1a}) into the dimensionless ITT equations (dropping the primes) which, from now on, will be the form used in this paper\,:
\beq{nv2}
\left(\partial_{t}+\bu\cdot\nabla\right)\bu + \nabla p &=& \alpha_{0}\bu + \Rn^{-1}\Delta\bu - \Rb\,\bu|\bu|^{2}\,,
\eeq
together with the incompressibility condition $\mbox{div}\,\bu = 0$. These operate on the unit periodic box $V_{d}= [0,\,1]^{d}$. The two Reynolds numbers $\Rn$ and $\Rb$ are defined as follows\,:
\beq{nv3}
\Rn &=& \frac{U_{0}L}{\nu}\,, \qquad\Rb = \frac{\beta U_{0}L}{\lambda^{2}}\,,
\label{eq:Re}
\eeq
taken together with the dimensionless frequency $\alpha_{0} = L\alpha U_{0}^{-1}$. The second choice of $U_{0}$ corresponds to $\Rn = 1$.


\section{Invariant scaling, time averages and length scales}\label{sec:scaling}

The incompressible NSEs possess the following well-known and powerful invariant scaling property involving an arbitrary parameter $\ell$:
\bel{inv}
\bx' = \ell^{-1}x ; \quad t' = \ell^{-2}t ;  \quad \bu = \ell^{-1}\bu' ;
\ee
which means that these equations are valid at every scale. The effect of this invariance is to scale the norms $\|\nabla^{n}\bu\|_{2m}$ defined by
\bel{norm1}
\|\nabla^{n}\bu\|_{2m} = \left(\int_{V_{d}}|\nabla^{n}\bu|^{2m}dV_{d}\right)^{1/2m}
\ee
in the following way:
\bel{norm2}
\|\nabla^{n}\bu\|_{2m} = \ell^{-1/\alpha_{m,n,d}}\|\nabla^{'n}\bu'\|_{2m} \,,
\ee
where $\alpha_{n,m,d}$ is defined by\footnote{The possible confusion caused by the labelling of the dimensionless frequency $\alpha_{0}$ and the exponents $\alpha_{n,m}$ is unfortunate, but we continue to use it to avoid the greater confusion of changing the notation from previous papers.} 
\bel{dim2a} 
\alpha_{n,m,d} = \frac{2m}{2m(n+1)-d}\,.
\ee
The $\alpha_{n,m,d}$ are a product of the invariance property (\ref{inv}). A dimensionless version of the norms defined in (\ref{norm2}) is given by
\bel{Fnmdef}
F_{n,m,d} := \nu^{-1}L^{1/\alpha_{n,m,d}}\|\nabla^{n}\bu\|_{2m}\,.
\ee
It has been shown that, for $d=2,\,3$, and for $n \geq 1$ and $1 \leq m \leq \infty$, weak solutions of the incompressible NSEs obey \cite{JDG2019,JDG2020}

\bel{dim1a}
\left<F_{n,m,d}^{(4-d)\alpha_{n,m,d}}\right>_{T} < \infty\,.
\ee
The angular brackets $\left<\cdot\right>_{T}$ represent the time average up to a time $T$, i.e., 
\bel{tadef}
\left<\cdot\right>_{T} = \frac{1}{T}\int_{0}^{T}\cdot\,d\tau\,.
\ee
We emphasize that these brackets represent a \textit{time average}, not a statistical average. When $n=0$ then $m$ is restricted by $3 < m \leq \infty$. An example familiar to the reader is the case $n=m=1$, in which case $(4-d)\alpha_{1,1,d} = 2$ with the cancellation of the factor of $4-d$. Then (\ref{dim1a}) yields the\footnote{The vorticity $\bom$ and the velocity gradient tensor $\nabla\bu$ are synonymous in $L^{2}$ when $\mbox{div}\,\bu=0$, but not in $L^{p}$ for $p>2$.} familiar bound on the time-averaged energy dissipation rate 
\bel{en2a}
\varepsilon = \nu L^{-d}\left<\int_{V}|\bom|^{2}dV_{d}\right>_{T} \leq \nu^{3}L^{-4}{\rm Re}^{3}.
\ee
With the inverse Kolmogorov length defined by $\lambda_{k}^{-4} = \varepsilon/\nu^{3}$, we obtain the conventional bound 
\bel{en2b}
L\lambda_{k}^{-1} \leq {\rm Re}^{3/4}\,.
\ee
Equation (\ref{dim1a}) thus expresses an infinite hierarchy of such bounds and can be looked upon as weighted space-time averages of all derivatives of the velocity field in every $L^{2m}$-norm. There is an informal analogy with the concept of wavelets\,: higher derivatives reflect the dynamics at small scales, while increasing the value of $m$ magnifies the larger amplitudes at each scale.

In \cite{JDG2019,JDG2020} it has also been shown how to define a set of inverse length scales associated with (\ref{dim1a}). Consider the set of $t$-dependent length-scales $\{\ell_{n,m,d}(t)\}$
defined by 
\bel{ls2}
\left(L\ell_{n,m,d}^{-1}\right)^{n+1} := F_{n,m,d}\,.
\ee
This definition takes into account the scaling of the domain volume $L^{d}$ which makes (\ref{ls2}) at the level of $n=m=1$ and $d=3$ consistent with the correct definition of the energy dissipation rate used to define the Kolmogorov length. Then we easily find that for Navier-Stokes weak solutions, when $n \geq 1$ and $1 \leq m \leq \infty$\,,
\bel{ls3}
\left<L\ell_{n,m,d}^{-1}\right>_{T} \leq c_{n,m,d}Re^{\frac{3}{(4-d)(n+1)\alpha_{n,m,d}}}\,.
\ee
When $d=3$ and $n=m=1$, then the exponent in (\ref{ls3}) is $\sthreequart$, which is consistent with estimates for the Kolmogorov length \cite{JDG2019,JDG2020}. Also, note that $(n+1)\alpha_{n,m,d} \to 1$ as $n,\,m \to\infty$\,.
\par\smallskip
Of course, it has been known for many years that solutions of the NSEs in $d=2$ dimensions are regular \cite{FMRT2001,DG1995}, but expressing (\ref{dim1a}) in integer $d$-dimensions $d=1,\,2,\,3$ rolls together into one line all the known two- and three-dimensional Navier-Stokes solution results, such as the class of weak solution $d=3$ time averages found by Foias, Guillop\'e and Temam \cite{FGT1981} in their pioneering paper in 1981.  It has been explained in Ref.~\cite{JDG2019} that, for a full existence and uniqueness proof in the $d=3$ case, a factor of $2\alpha_{n,m,3}$ would be needed in the exponent in (\ref{dim1a}). However, no evidence exists for the existence of bounds with this necessary factor of $2$. It is possible that the Leray-Hopf weak solutions are all that exist. 

By inspection it is clear that the ITT equations respect the invariant scaling possessed by the NSEs equations, apart from the linear-pumping term. However, there is a significant literature on a more general class of equations where the $\bu|\bu|^{2}$-term is replaced by $\bu|\bu|^{2\delta}$, which is the case in the Brinkman-Forchheimer extended Darcy model arising in porous media. The paper by Titi and Trabelsi \cite{TT2022} contains a wide literature survey, but we also refer the reader to \cite{Anton2010,Cai2008,Celebi2006,Liu2007,LSST2015,LST2015,MTT2016,Payne1999,BS2008,Ugurlu2008,Wang2008,You2012}. When $\delta > 1$\,, the invariant scaling property of the NSEs is broken. This leads to the bounding of time-averaged norms, higher than those available to ITT, which eventually lead to the regularity of solutions in the $d=3$ case. The work in our paper is different in two important respects. Firstly, we focus on the critical value $\delta = 1$ and thus remain faithful to the scaling property in (\ref{inv}). Secondly, the ITT equations in (\ref{nv2}) have a linear term $\alpha_{0}\bu$ which, while trivial in a purely functional setting, is nevertheless physically important in the creation of equilibrated or frozen states, which appear to dominate the dynamics in our DNSs. 

The parallel scaling properties of the ITT equations and the NSEs suggest that the exponents $\alpha_{n,m,d}$ in (\ref{dim2a}) should be the same in both cases. Therefore, taking into account the factor of $4-d$ in the exponent in the $d=2$ case, we define 
\bel{Pnmdef}
P_{n,m} = \|\nabla^{n}\bu\|_{2m}^{2\alpha_{n,m,2}} \,, \qquad(d=2)\,,
\ee
where 
\bel{alpha2def}
\alpha_{n,m,2} = \frac{m}{m(n+1)-1}
\ee
and a set of inverse length scales equivalent to (\ref{ls3}) 
\bel{ls4}
\left<L\ell_{n,m,2}^{-1}\right>_{T} \leq c_{n,m,2}\left<P_{n,m}\right>_{T}^{\frac{1}{2(n+1)\alpha_{n,m,2}}}\,.
\ee
When $d=3$ 
\bel{Qnmdef}
Q_{n,m} = \|\nabla^{n}\bu\|_{2m}^{\alpha_{n,m,3}}\,, \qquad(d=3)\,,
\ee
where 
\bel{alpha3def}
\alpha_{n,m,3} = \frac{2m}{2m(n+1)-3}
\ee
and a set of inverse length scales equivalent to (\ref{ls3}) 
\bel{ls5}
\left<L\ell_{n,m,3}^{-1}\right>_{T} \leq c_{n,m,3}\left<Q_{n,m}\right>_{T}^{\frac{1}{(n+1)\alpha_{n,m,3}}}\,.
\ee
General bounds on $\left<P_{n,m}\right>_{T}$, expressed as a function of $\alpha_{0}$, $\Rn$ and $\Ra$, are estimated in \S\ref{sec:sumd2}. For $d=3$, a narrower class of bounds on $\left<Q_{n,m}\right>_{T}$ has been given in \S\ref{sec:sumd3}.


\section{Energy estimates}\label{sec:energy}

In keeping with standard Navier-Stokes notation, we define $n$ derivatives of $\bu$ in $L^{2}(V_{d})$ as 
\bel{Hndef}
H_{n} = \int_{V_{d}}|\nabla^{n}\bu|^{2}dV_{d}\,.
\ee
Given the close relationship between the ITT equations and the incompressible NSEs, a formal approach is taken following the method used by \cite{FGT1981}. In this notation the energy $H_{0}$ and and enstrophy $H_{1}$ are
\beq{enenstA}
H_{0} &=& \int_{V_{d}}|\bu|^{2}dV_{d}\,; \\
H_{1} &=& \int_{V_{d}}|\nabla\bu|^{2}dV_{d} = \int_{V_{d}}|\bom|^{2}dV_{d}\,. \label{enenstB}
\eeq
Using standard methods \cite{Leray1934,FMRT2001}, a Leray-type energy inequality is easily derived 
\bel{ee1}
\shalf \dot{H}_{0} + \Rn^{-1}H_{1} + \Rb \int_{V}|\bu|^{4}dV_{d} \leq \alpha_{0}H_{0} 
\ee
from which we drop the $H_{1}$-term\footnote{Poincar\'e's inequality cannot be applied because, unlike the NSEs, the spatial average of $\bu$ is not zero.} and apply a H\"older inequality to the $L^{4}$-term to produce a simple differential inequality for $H_{0}$ 
\bel{ee2a}
\shalf \dot{H}_{0} \leq \alpha_{0}H_{0} - \Rb H_{0}^{2}\,.
\ee
Thus, equilibration of the right hand side occurs at 
\bel{ee2b}
H_{0,\,equil} = \alpha_{0}\Rb^{-1}\,:=\,\Ra \,,
\ee
where we designate $\Ra$ as the activity parameter. By using the time-average definition in (\ref{tadef}), from (\ref{ee1}) we find the following estimates
\beq{ee3a}
\left<H_{0}\right>_{T} &\leq& \Ra\,, \\
\left<H_{1}\right>_{T} &\leq& \alpha_{0}\Ra\Rn\,,\label{ee3b}
\eeq
{together with the time average of both the $L^{4}$-norm and the ratio of $H_{1}$ to $H_{0}$
\bel{ee3c}
\left<\int_{V}|\bu|^{4}dV_{d}\right>_{T}  \leq \Ra^{2}\,,\qquad\qquad 
\left<\frac{H_{1}}{H_{0}}\right>_{T} \leq \alpha_{0}\Rn\,.
\ee}
The inequalities (\ref{ee3a})--(\ref{ee3c}) each have an $O\left(T^{-1}\right)$ correction term that will be dropped from now on. These results are true in every dimension. With respect to initial data on $H_{0}$\,: 
\ben\itemsep 1mm
\item For initial data $H_{0}(0) > \Ra$\,, the sign of $\dot{H}_{0}$ in this region is negative and thus $H_{0}$ decreases to $\Ra$. 

\item For initial data $H_{0}(0) < \Ra$\,, {$\dot{H}_{0}$ can take either sign. If $\dot{H}_{0}(t) > 0$ then $H_{0}$ will grow to reach $\Ra$ (a frozen state) but cannot pass through it. If $\dot{H}_{0}(t) < 0$ then $H_{0}$ will decay.}
\een
The choice of $n=m=1$ makes $(4-d)\alpha_{1,1,d} = 2$, whereupon the factor of $4-d$ cancels. Thus
\bel{ee4a}
\left<P_{1,1}\right>_{T} \leq \alpha_{0}\Ra\Rn\,,
\ee
and 
\bel{ee4b}
\left<Q_{1,1}\right>_{T} \leq \alpha_{0}\Ra\Rn \,.
\ee


\section{Numerical Methods}\label{sec:num}
\begin{table}[!h]
		\centering
		\begin{tabular}{|c| c| c| c|c|c|c| c|}
			\hline
			Run & $d$ &  $N$& $\delta t$&$\nu$&$\alpha$&$\beta$ 
			\\
			\hline
			A1& 2& 2048&$2\times 10^{-4}$& $2.87\times10^{-1}$&$10\times10^{1}$&$5$\\
			\hline
			A2& 2& 2048&$2\times 10^{-4}$& $1.41\times10^{-1}$&$10\times10^{1}$&$5$\\
			\hline
			A3& 2& 2048&$2\times 10^{-4}$& $7.07\times10^{-2}$&$10\times10^{1}$&$5$\\
			\hline
			A4& 2& 2048&$2\times 10^{-4}$& $3.53\times10^{-2}$&$10\times10^{1}$&$5$\\
			\hline
			A5& 2& 2048&$2\times 10^{-4}$& $2.36\times10^{-2}$&$10\times10^{1}$&$5$\\
			\hline
			A6& 2& 2048&$2\times 10^{-4}$& $1.77\times10^{-2}$&$10\times10^{1}$&$5$\\
			\hline
			A7& 2& 2048&$2\times 10^{-4}$& $1\times10^{-2}$&$10\times10^{1}$&$5$\\
			\hline
			A8& 2& 2048&$2\times 10^{-4}$& $8.8\times10^{-2}$&$10\times10^{1}$&$5$\\
			\hline
			F1& 2& 2048&$2\times 10^{-4}$& $6.2\times10^{-1}$&$1$&$1$\\
			\hline
			F2& 2& 2048&$2\times 10^{-4}$& $1.2\times10^{-1}$&$1$&$1$\\
			\hline
			F3& 2& 2048&$2\times 10^{-4}$& $6.0\times10^{-2}$&$1$&$1$\\
			\hline
			F4& 2& 2048&$2\times 10^{-4}$& $3.0\times10^{-2}$&$1$&$1$\\
			\hline
			F5& 2& 2048&$2\times 10^{-4}$& $1.5\times10^{-2}$&$1$&$1$\\
			\hline
			F6& 2& 2048&$2\times 10^{-4}$& $7.0\times10^{-3}$&$1$&$1$\\
			\hline
			F7& 2& 2048&$2\times 10^{-4}$& $3.1\times10^{-3}$&$1$&$1$\\
			\hline
			B1& 3&  512&$1\times10^{-3}$ &$5\times10^{-1}$&$1\times10^{1}$&$1\times10^{-1}$\\
			\hline
			B2& 3&  512&$1\times10^{-3}$ &$5\times10^{-2}$&$1\times10^{1}$&$1\times10^{-1}$\\
			\hline
			B3& 3&  512&$5\times10^{-4}$ &$1\times10^{-2}$&$1\times10^{1}$&$1\times10^{-1}$\\
			\hline
		\end{tabular}
		\caption{\label{tab:widgets} The parameters for our DNSs: $d$ is the dimension, $N^d$ the number of collocation points, and $\delta$ the time step.For all our runs, $\lambda=1$. Given these parameters, the Reynolds numbers follow from Eqs.~(\ref{a1b}) and	(\ref{eq:Re}). Parameters for other runs are given in the Supplemental Material.}
		\label{tab:parameters}
	\end{table}
For our DNS of the $d$-dimensional ITT Eq.~(\ref{a1a}), we use a Fourier pseudospectral method \cite{canuto2012spectral} on periodic domains (a square in $d=2$ and a cube in $d=3$), with sides of length $L=2\pi$, and $N^{d}$ collocation points. We employ the second-order exponential time-difference scheme, ETDRK2, for time evolution in Fourier-space~\cite{cox2002}. We list the parameters for our DNS runs in Table \ref{tab:parameters}.  Our simulations are numerically robust insofar as the Courant-Friedrichs-Lewy number, $C<C_{max}$, where $C=\left(\sum_{i=1}^{i=d}U_{i}/h_{i}\right)\delta t$ with $U_{i}=~\sup_{x}u_{i}(\bx,t)$ 
and $h_{i}$ the minimum grid spacing between two collocation points. Typically $C_{max}=1$ for our explicit numerical schemes and for all our simulations (Table~\ref{tab:parameters}) $C \simeq  0.62\pm0.04$.

 The dimensional version of the ITT equations (\ref{a1a}) has four parameters $\lambda,\,\alpha,\,\nu$ and $\beta$ which reduce to the three dimensionless numbers $Re_{\nu},\,Re_{\beta}$ and $\alpha_{0}$ in the non-dimensionalized version (\ref{nv2}). As explained in (\ref{a1b}), we have found it convenient to define the typical velocity field $U_{0}$ in two particular ways\,: $U_{0}= \sqrt{\alpha/\beta}$ and $U_{0} = \nu/L$. The latter case restricts $Re_{\nu}$ to the value $Re_{\nu}=1$ but allows us to explore a more diverse range of $\alpha_{0}$ and $Re_{\beta}$. 

\section{Summary of results in the $d=2$ case}\label{sec:sumd2}

The methods used in the analysis sections of this paper are based on the differential  inequalities explained in  \ref{appinequal}. The proof of the results in the following subsections are given in \ref{appA}. Within these estimates, various multiplicative constants $c,\,c_{m}$ and $c_{n,m}$ appear, which should be read as generic constants that may differ from line to line. These constants are algebraic in $n,m$ but are not usually given explicitly\,: see \ref{appinequal}. {We remark that none of the bounds displayed in the following sections are saturated, although more drastic initial conditions might get closer to this state.}
{For simplicity, it is also assumed that $\Rn \gg \Rb$. In the choice $U_{0} = \nu/L$ where $\Rn=1$, in which case estimates can be re-calculated from the material in the Appendices.}

\subsection{Estimates for $\left<P_{n,m}\right>_{T}$}\label{subsec:est2}

By using the definition of $P_{n,m}$ in (\ref{Pnmdef}) as our guide, for $m=1$, we have $\alpha_{n,1,2} = 1/n$\,, so
\bel{P1a}
\left<P_{n,1}\right>_{T} = \left<H_{n}^{1/n}\right>_{T} \,.
\ee
The estimate for $\left<P_{1,1}\right>_{T}$ in (\ref{ee4a}) can be used to compute a series of other inequalities. The proofs can be found in \ref{appA}. Below is a summary\,:
\par\smallskip\noindent
(i) Firstly, we wish to estimate $\left<P_{n,1}\right>_{T}$ for $n \geq 1$. In inequality (\ref{P9}) in \ref{P2sect}, for $n=2$ it is shown that
\bel{P1c}
\left<P_{2,1}\right>_{T} \leq c\,\alpha_{0}\left(\alpha_{0}\Ra\Rn^{3}\right)^{1/2} \,.
\ee
\par\smallskip\noindent
(ii) More generally, in inequality (\ref{Pn6}) \ref{Pnsect}, it is shown that for $n\geq 2$
\beq{P2a}
\left<P_{n,1}\right>_{T} &\leq& c_{n,1}\alpha_{0}^{2/n}\left(\alpha_{0}\Ra\Rn^{3}
\right)^{\frac{n-1}{n}} \,.
\eeq
\par\smallskip\noindent
(iii) In inequality (\ref{P1mb}) in \ref{P1m} it is shown that
\bel{P2b}
\left<P_{1,m}\right>_{T} = \left<\|\nabla\bu\|_{2m}^{\frac{2m}{2m-1}}\right>_{T}\,,
\ee
displayed in row 3 of Fig.~\ref{fig:fig1}, satisfies
\beq{P2c}
\left<P_{1,m}\right>_{T}  
&\leq& c_{m}\left(\alpha_{0}\Rn\right)^{\frac{3m-2}{2m-1}}\Ra^{\frac{m}{2m-1}}\,.
\eeq
{(\ref{P1mc}) shows that in the limit $m\to\infty$, we also have 
\bel{P2d1}
\left<\|\nabla\bu\|_{\infty}\right>_{T} \leq c\, \left(\alpha_{0}\Rn\right)^{3/2}\Ra^{1/2}\,.
\ee}
\par\smallskip\noindent
(iv) With the definition 
\bel{P2d2}
\left<P_{0,m}\right>_{T} = \left<\|\bu\|_{2m}^{\frac{2m}{m-1}}\right>_{T}\,,
\ee
 from (\ref{P12}) in \ref{uinf}, it is shown that, for $m > 2$\,, 
\bel{P2e}
\left<P_{0,m}\right>_{T} \leq c\,\Ra^{\frac{m}{m-1}}\left(\alpha_{0}\Rn\right)^{\frac{m-2}{m-1}}\,.
\ee
(\ref{P13}) shows that, in the limit $m\to \infty$, we have 
\bel{P2f}
\left<P_{0,\infty}\right>_{T} = \left<\|\bu\|_{\infty}^{2}\right>_{T} \leq c\,\alpha_{0}\Ra\Rn \,.
\ee
\par\smallskip\noindent
(v) In inequality (\ref{Pnms5}) in \ref{Pnmsect} it is shown that, for $n \geq 2$,
\bel{P2g}
\left<P_{n,m}\right>_{T} \leq c_{n,m}\alpha_{0}^{\frac{2m}{m(n+1)-1}}
\left(\alpha_{0}\Ra\Rn^{3}\right)^{\frac{mn-1}{m(n+1)-1}}\,.
\ee
\par\smallskip\noindent

{
\subsection{Regularity\,: exponential bounds in $d=2$ dimensions}\label{subsec:reg1}

The difference between the $2d$ and $3d$ NSEs lies in the absence of the vortex stretching term $\bom\cdot\nabla\bu$, which is zero when $d=2$. For the ITT equations (\ref{nv2}), the evolution equation for the vorticity is
\bel{reg1}
\left(\partial_{t} + \bu\cdot\nabla\right)\bom = \alpha_{0}\bom + \Rn^{-1}\Delta\bom - \Rb\,\mbox{curl}\left(\bu|\bu|^{2}\right)
\ee
in which the $\mbox{curl}\left(\bu|\bu|^{2}\right)$-term appears to create another form of vortex stretching. Specifically we have
\beq{reg2}
\shalf \dot{H}_{1} &=& \alpha_{0}H_{1} - \Rn^{-1}H_{2} - \Rb\I \bom\cdot \mbox{curl}(\bu|\bu|^{2})dV\nonumber\\
&\leq& \alpha_{0}H_{1} - \shalf \Rn^{-1}H_{2} + \shalf \Rb^{2}\Rn\I |\bu|^{6}dV\,,
\eeq
where we have integrated by parts and have then used a H\"older inequality. Then we use a Gagliardo-Nirenberg inequality
\bel{reg3}
\|\bu\|_{6} \leq c\, \|\nabla\bu\|_{2}^{\sthird}\|\bu\|_{4}^{\twothirds}\,,
\ee
which can be re-expressed as
\bel{reg4}
\I |\bu|^{6}dV \leq c\,H_{1}\I |\bu|^{4}dV\,.
\ee
Thus we can write
\bel{reg5}
\shalf \dot{H}_{1} \leq \left(\alpha_{0} + \shalf \Rb^{2}\Rn \I|\bu|^{4}\,dV\right)H_{1}\,,
\ee
whence we can use the extra piece of information afforded to us by the bounded time integral expressed in (\ref{ee3c}). Thus (\ref{reg5}) can be integrated to become 
\beq{reg7}
H_{1}(T) &\leq& H_{1}(0) \exp \left\{\int_{0}^{T}\left(\alpha_{0} + c\,\Rb^{2}\Rn\|\bu\|_{4}^{4}\right)\,d\tau\right\}\nonumber\\
&\leq& H_{1}(0) \exp\left\{\alpha_{0}\left(1 + c\,\Rb^{2}\Rn\Ra^{2}\right)T\right\}\,,
\eeq
which is finite for every finite $T$. Control over the $H_{1}$-norm establishes global regularity in this $2d$ case but not a global attractor, which requires a uniform bound for all $t$. 
}

\subsection{Spectral energy budget}\label{subsec:seb}

\begin{figure}[!h]
	\resizebox{\linewidth}{!}{
 	\includegraphics[scale=1]{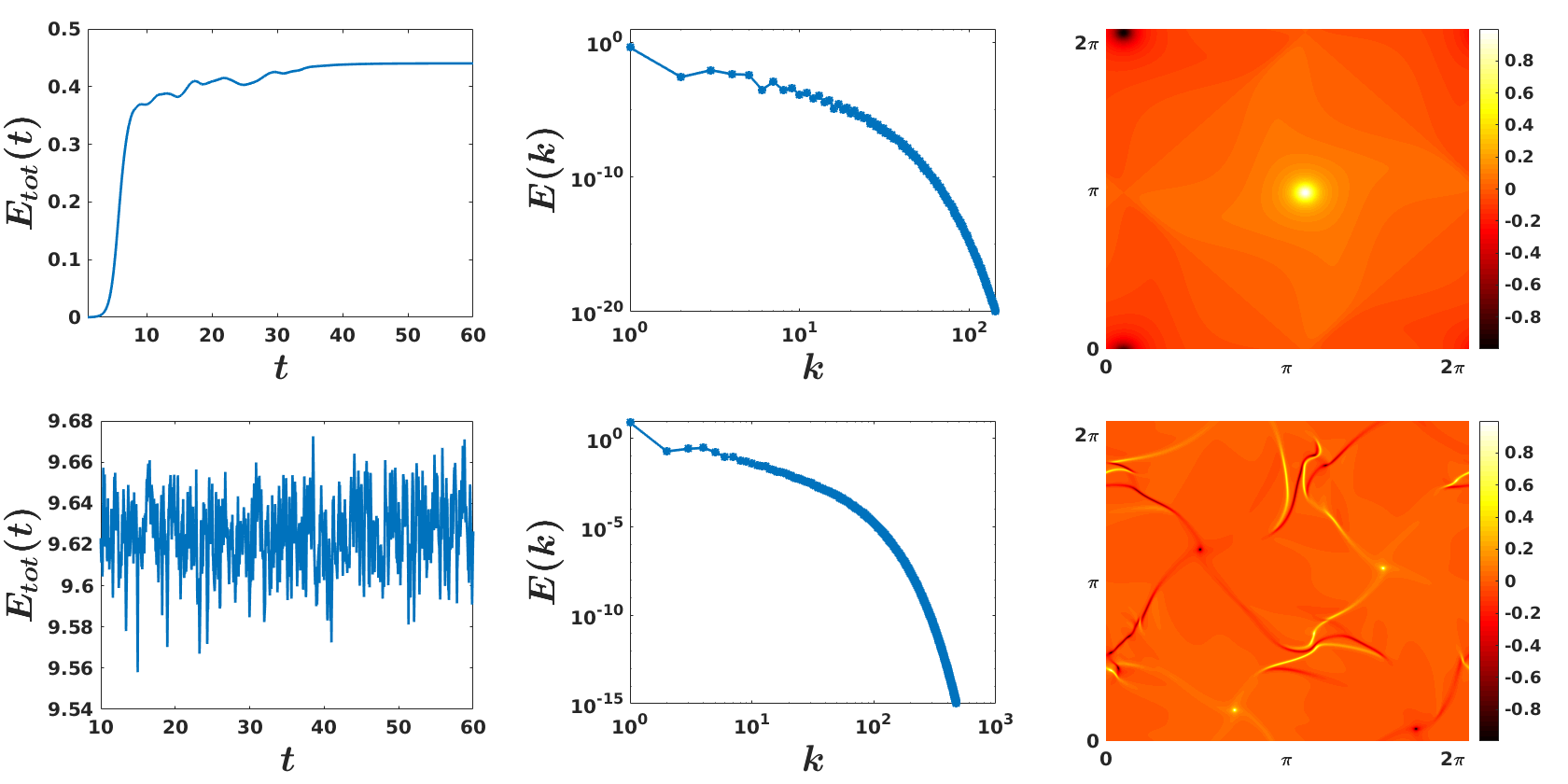}
 	\put(-380,200){\rm {\bf(a)}}
 	\put(-240,200){\rm {\bf(b)}}
 	\put(-120,200){\rm {\bf(c)}}
 	\put(-380,95){\rm {\bf(d)}}
 	\put(-240,95){\rm {\bf(e)}}
 	\put(-120,95){\rm {\bf(f)}}
 	}
	\caption{{Plots for runs F7 (row 1) and A6 (row2): column (1) contains plots versus the time $t$ of the total energy $E_{tot}(t)$; column (2) contains log-log plots versus $k$ of the energy spectrum $E(k)$; the  filled contour plots in column (3) are of $\bom$ at a representative time.}}
	\label{fig:2d_spectra}
	\end{figure}

\begin{figure}[!h]
    \resizebox{\linewidth}{!}{
    \includegraphics[scale=1]{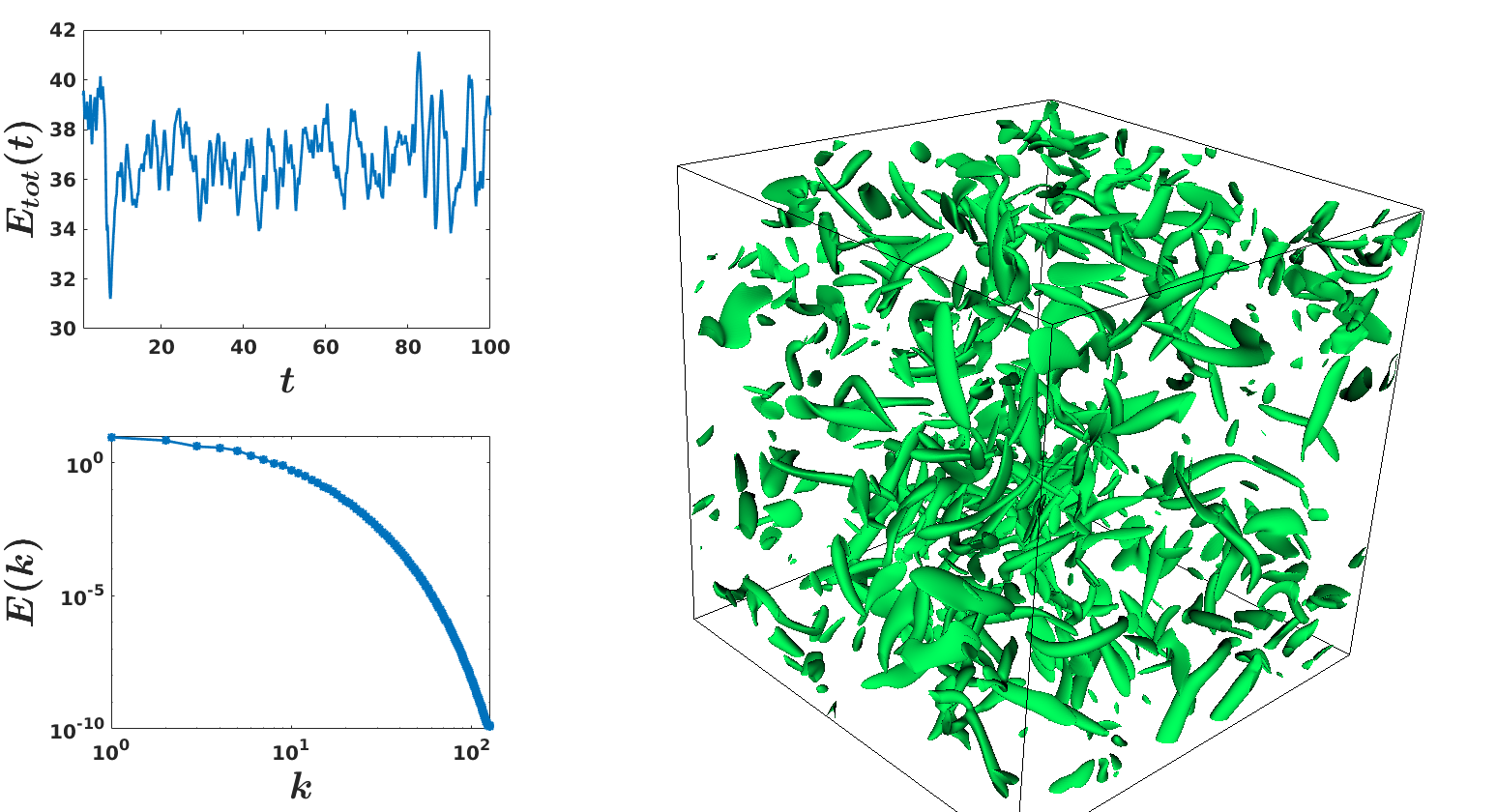}
    \put(-380,200){\rm {\bf(a)}}
    \put(-380,100){\rm {\bf(b)}}
    \put(-200,200){\rm {\bf(c)}}
    }
    \caption{{(Colour online) For run B2;  plot versus $t$ of the total kinetic energy $E_{tot}(t)$ (column 1, row 1); log-log plot versus $k$ of the energy spectrum $E(k)$ (column 1, row2); iso-surfaces of the modulus of the vorticity field (column 2). }}
    \label{fig:3d_spectra}
\end{figure}

\begin{figure*}[h!]
    	\resizebox{\linewidth}{!}{
    		\includegraphics[scale=1]{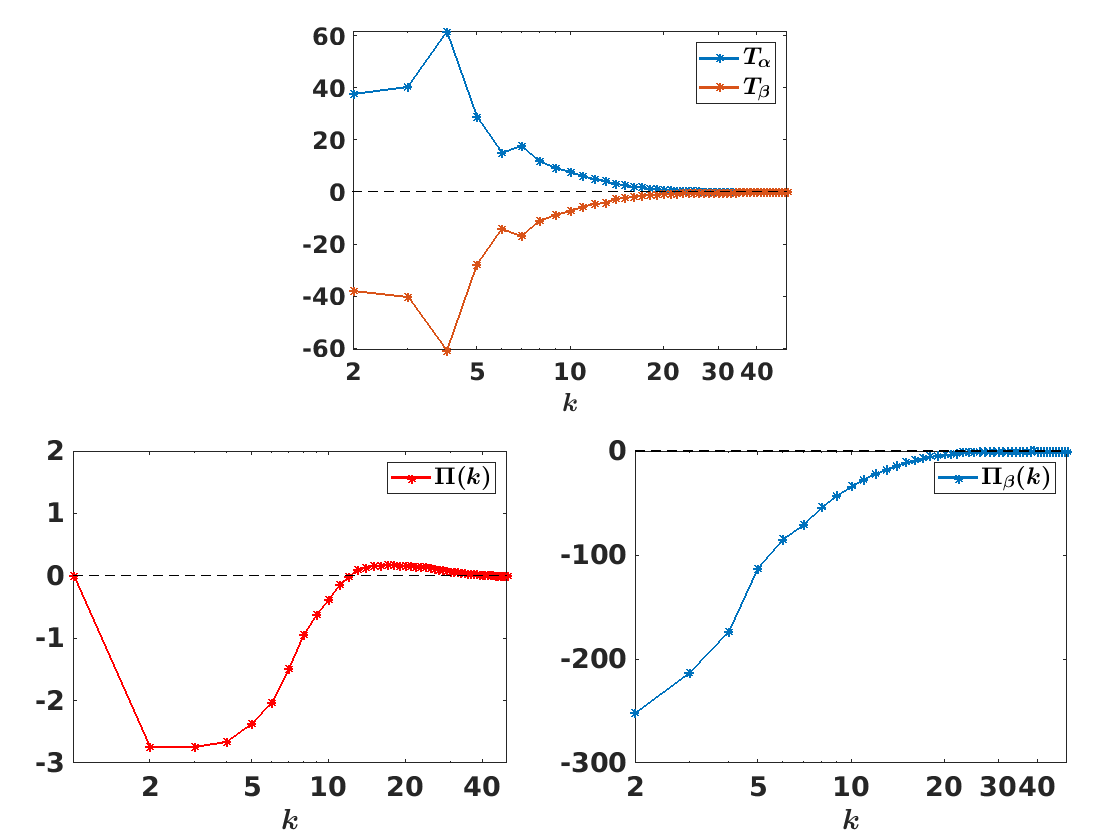}
    		\put(-180,200){\rm {\bf(a)}}
    		\put(-240,100){\rm {\bf(b)}}
    		\put(-100,100){\rm {\bf(c)}}
    		    		}
    	\caption{({Color online) Semi-log plots versus $k$ of (a) $T_{\alpha}(k)$ (blue) and $T_{\beta}(k)$ (red), (b) $\Pi(k)$ (red), and (c) $\Pi_{\beta}(k)$ (blue), for an illustrative simulation (A4) from Table. \ref{tab:parameters}}}
		\label{fig:spectra_budget}
    \end{figure*}

{We discuss the role of contributions to the energy from the nonlinear terms in the ITT equations. The shell-averaged energy spectrum is defined by
   \begin{equation}\label{eq:spectra}
       \mathcal{E}(k)=\frac{1}{2} \sum_{k'=k-1/2}^{k'=k+1/2}\left<\widetilde{\bu}(\bk')\cdot\widetilde{\bu}(-\bk')\right>_{t}\,,
   \end{equation}
   where a tilde denotes a spatial Fourier transform and the wave vectors $\bk$ and $\bk'$ have 
   moduli $k$ and $k'$, respectively. In Figs.~\ref{fig:2d_spectra} and \ref{fig:3d_spectra} we give illustrative plots of the time series of the total kinetic energy
   \bel{Edef}
   E_{tot}(t)\equiv \shalf L^{-d} \int_{V_{d}}|\bu|^{2}\,dV_{d}
   \ee
   together with $\mathcal{E}(k)$, filled contour plots of the vorticity (in $d=2$) and isosurfaces of the modulus of the vorticity (in $d=3$). From Fig.~\ref{fig:2d_spectra}(a) we surmise the existence of temporally frozen states in $d=2$, because $dE_{tot}/dt \simeq 0$ after the period of initial transients. Both frozen and turbulent states can also be characterized by the filled contour plots of the vorticity presented in column 3 of Fig.~\ref{fig:2d_spectra}; these states depend on $\nu$, $\alpha$, and $\beta$, but not on the initial conditions\,; a change in the initial conditions changes the time for which initial transients last. We also find statistically steady states, with fluctuations in the total energy, in both $d=2$  (Fig.~\ref{fig:2d_spectra}(d)) and $d=3$ (Fig.~\ref{fig:3d_spectra}(a)), similar to those observed in homogeneous isotropic fluid turbulence.
   \par\smallskip
   The $k$-dependent energy budget is given by
   \begin{equation}{\label{eq:enbudget}}
   	\partial_{t}\mathcal{E}(k)=T(k)-T_{\beta}(k)+T_{\alpha}(k) -2\Gamma_{0}k^{2}\mathcal{E}(k)
   \end{equation}
   where $T_{\alpha}(k)$, the spectral energy contribution from the $\alpha$ term, and $T(k)$ and $T_{\beta}(k)$ are, respectively, the $k$-shell averaged contributions from the advective and cubic terms in the ITT equations, are:
   	\begin{align}{\label{eq:enonlin}}
   		T(k)&=-\lambda\sum_{k'=k-1/2}^{k'=k+1/2}\left<\widetilde{\bu}(-\bk)\cdot\textbf{P}\left(\bk'\right)
   		\cdot\widetilde{\left(\bu\cdot\nabla\bu\right)}(\bk')\right>_{t}\,; \\ 
   	   T_{\beta}(k)&=\beta\sum_{k'=k-1/2}^{k'=k+1/2}\left< \widetilde{\bu}(-\bk)\cdot\textbf{P}(\bk')\cdot\widetilde{(\vert\bu\vert^{2}\bu)}(\bk')\right>_{t} \,;\\
   	   T_{\alpha}&=-2\alpha\mathcal{E}(k) \,;
   	\end{align}
    here, 
    \bel{Pdef}
    P_{ij}(\bk)=\delta_{ij}-\frac{k_{i}k_{j}}{k^{2}}
    \ee
    is the transverse projector and $\left< \cdot \right>_{t}$ the average over time $t$. The flux of energy arising from the advective and the cubic terms are, respectively:
    \begin{align}{\label{eq:flux}}
    \Pi(k)&=-\sum_{k'=0}^{k'=k}T(k')\,; \\
    \Pi_{\beta}(k)&=-\sum_{k'=0}^{k'=k}T_{\beta}(k')\,.   
    \end{align}
    In Fig.~\ref{fig:spectra_budget}(a), we plot versus $k$, $T_{\beta}(k)$ and $T_{\alpha}(k)$; we omit the points $T_{\alpha}(1)$ and $T_{\beta}(1)$ as $T_{\alpha}(1)>>T_{\alpha}(k)$ and $T_{\beta}(1)<<T_{\beta}(k)$. We note from $T_{\beta}(k)$ that the cubic term is an energy sink and is dominant at small-$k$ modes; by contrast, the $\alpha$ term acts as an energy source.
   \par\smallskip
   In Figs.~\ref{fig:spectra_budget}(b) and (c) we give plots versus $k$ of the fluxes $\Pi(k)$ and $\Pi_{\beta}(k)$, respectively. The flux associated with the cubic term $\Pi_{\beta}(k)<0$, as we expect from an energy sink. However, the advective term is neither a sink nor a source of energy. Therefore, the total area under the curve in Fig.~\ref{fig:spectra_budget} (a) is zero.  We also note that there is no region of $k$ for which $\Pi(k)\simeq$ a constant. This indicates the absence of a conventional inertial range.}
    
\subsection{Numerical results for $P_{n,m}$}\label{subsec:num2}

\begin{figure*}[!]


		
		\begin{tikzpicture}
		\node[anchor=south west,inner sep=0] at (0,0)
		{\includegraphics[width=0.5\linewidth]{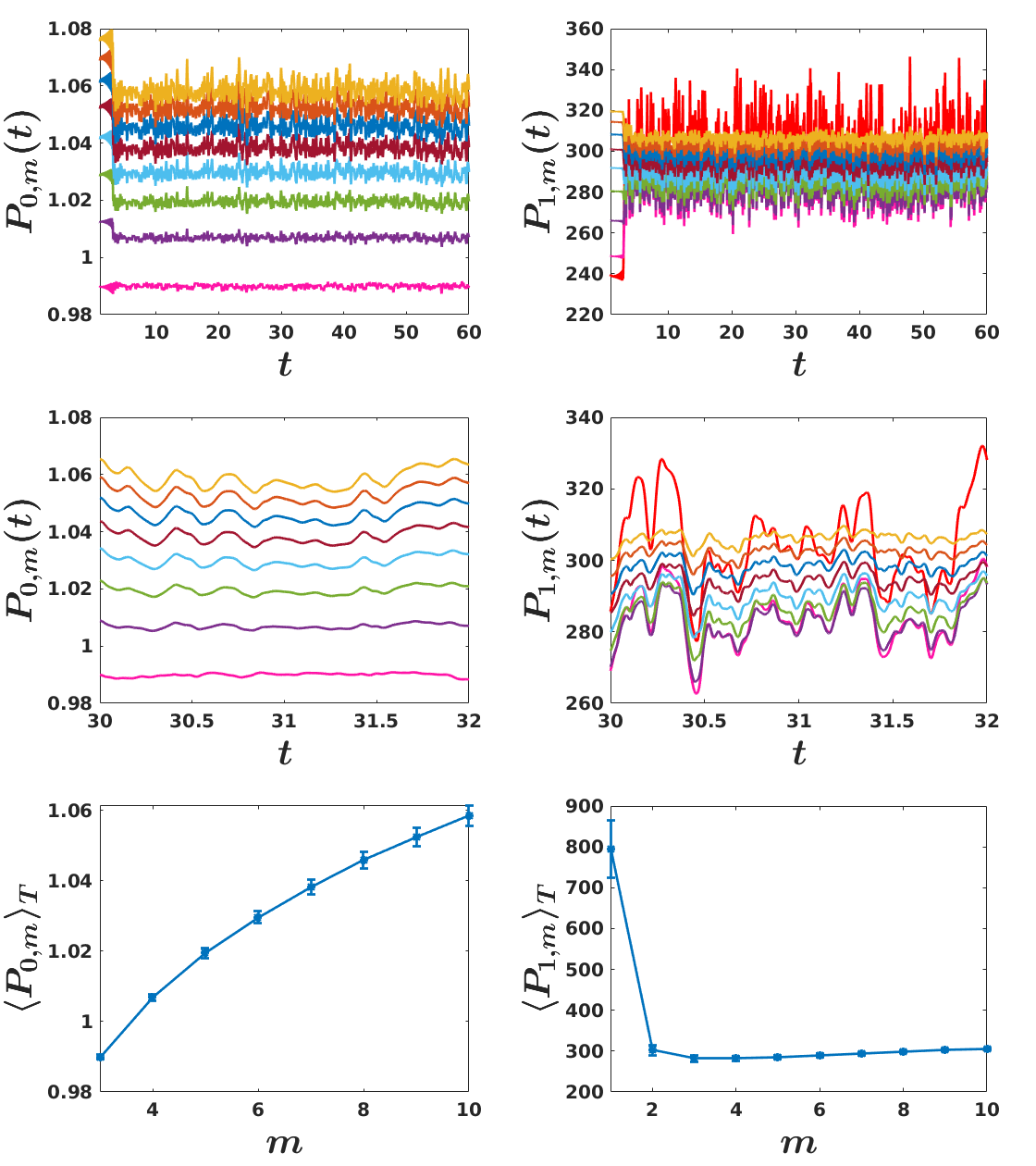}
		\put(-100,225){\large A}
		\put(90,225){\large B}
	
		\includegraphics[width=0.5\linewidth]{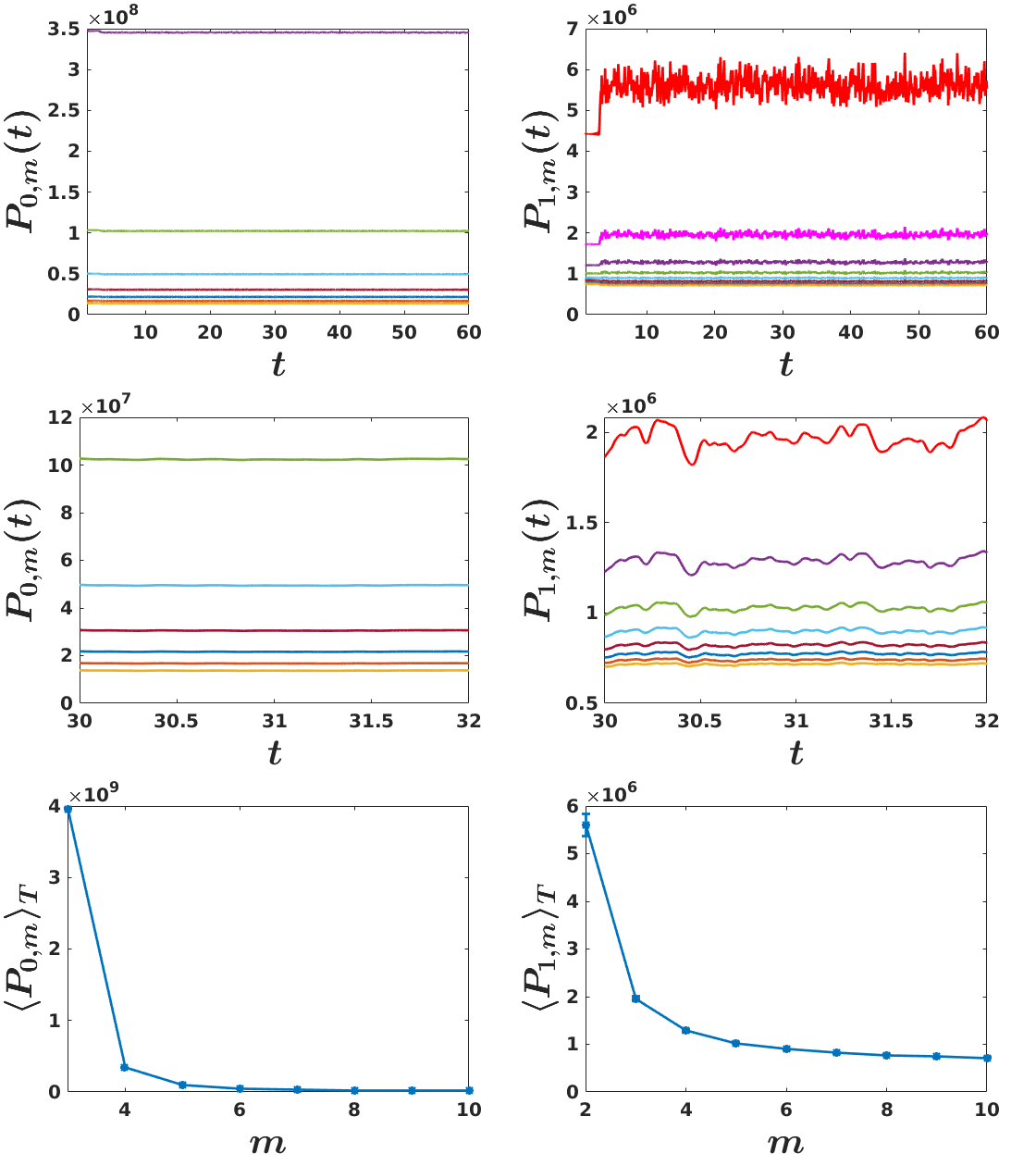}};
		\draw[line width=0.2mm,black, -] (0,0.0) -- (6.7,0.0);
		\draw[line width=0.2mm,black, -] (0,0.0) -- (0.0,7.8);
		\draw[line width=0.2mm,black, -] (6.7,0.0) -- (6.7,7.8);
		\draw[line width=0.2mm,black, -] (0,7.8) -- (6.7,7.8);
		\draw[line width=0.1mm,blue, dashed] (1.825,5.55) -- (0.7,4.9);
		\draw[line width=0.1mm,blue, dashed] (2.0,5.55) -- (3.1,4.9);
		\draw[line width=0.1mm,blue, dashed] (5.2,5.55) -- (4.0,4.9);
		\draw[line width=0.1mm,blue, dashed] (5.3,5.55) -- (6.4,4.9);
		\draw[line width=0.2mm,black, -] (6.75,0) -- (13.6,0);
		\draw[line width=0.2mm,black, -] (6.75,0) -- (6.75,7.8);
		\draw[line width=0.2mm,black, -] (6.75,7.8) --(13.6,7.8);
		\draw[line width=0.2mm,black, -] (13.6,0) --(13.6,7.8);
		\draw[line width=0.1mm,blue, dashed] (8.625,5.55) -- (7.5,4.9);
		\draw[line width=0.1mm,blue, dashed] (8.8,5.55) -- (9.9,4.9);
		\draw[line width=0.1mm,blue, dashed] (11.9,5.55) -- (10.8,4.9);
		\draw[line width=0.1mm,blue, dashed] (12.1,5.55) -- (13.2,4.9);
		\end{tikzpicture}
		\caption{(colour online) Illustrative plots for $U_{0} = \sqrt{\alpha/\beta}$ (panel A) and $U_{0} = \nu/L$ (panel B) for run A6 (see Table \ref{tab:parameters}): First and second rows\,: plots versus $t$ of $P_{0,m}$ and $P_{1,m}$\,; the plots in the second row are expanded versions of small segments of the plots in the first row. Third row\,: Plots versus $m$ of $\left< P_{0,m} \right>_{T}$ and $\left< P_{1,m} \right>_{T}$. Curves for $m=2, 3, 4, 5, 6, 7, 8, 9$, and $10$ are drawn in red, pink; violet, green, cyan, maroon, blue, orange, and yellow, respectively. Similar plots for other representative runs are given in the Supplemental Material.}\label{fig:fig1}
		\end{figure*}
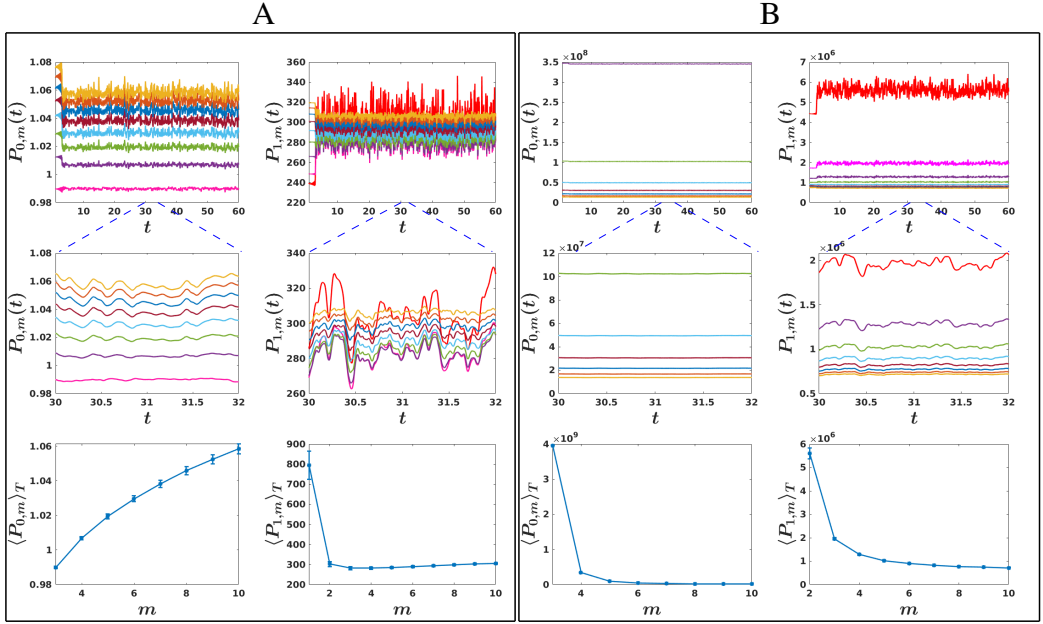

We now present plots in Fig.~\ref{fig:fig1} of $P_{n,m}(t)$ versus time $t$ for the two values $n=0$ and $n=1$, with a sequence of values of $m =  1,\, \ldots,\, 10$. H\"older's inequality insists that, for fixed $n$, the norms $\|\cdot\|_{2m}$ must be ordered with increasing $m$, such that $\|\cdot\|_{2m} \leq \|\cdot\|_{2(m+1)}$; but the $\alpha_{n,m,d=2}$ \textit{decrease} as $m$ increases. Thus, it is technically possible for the $P_{n,m}$ to be ordered either way\,: i.e.,  an increasing regime $P_{n,m} \leq P_{n,m+1}$ or a decreasing regime $P_{n,m} \geq P_{n,m+1}$. The latter regime was originally observed numerically for the $d=3$ NSEs~\cite{DKGGPV2013,GDKGPV2014} and was also discussed in \cite {GGKPPPPSS2016}, although no obvious reason for this particular ordering was deduced. Moreover, no crossing of curves that represented different values of $m$ has been observed. 
\par\smallskip
{For the $d=2$ case,} when $U_{0} = \sqrt{\alpha/\beta}$ (panel A), we have observed both regimes, but only the decreasing regime $P_{n,m} \geq P_{n,m+1}$ when $U_{0} = \nu/L$ (panel B). {This is illustrated in Fig.~\ref{fig:fig1} for run A6\,:} In panel A ($U_{0} = \sqrt{\alpha/\beta}$) the plots of $P_{1,m}$ can cross each other at different times, as we can see clearly in the expanded plots in the second rows\,; such crossings do not occur in panel B ($U_{0} = \nu/L$). Furthermore, the plots versus $m$ of $\left< P_{1,m} \right>_{T}$ (third row) decreases monotonically with increasing $m$ in panel B but not in panel A.
\par
\begin{figure*}[!]


		
		\begin{tikzpicture}
		\node[anchor=south west,inner sep=0] at (0,0)
		{\includegraphics[width=0.5\linewidth]{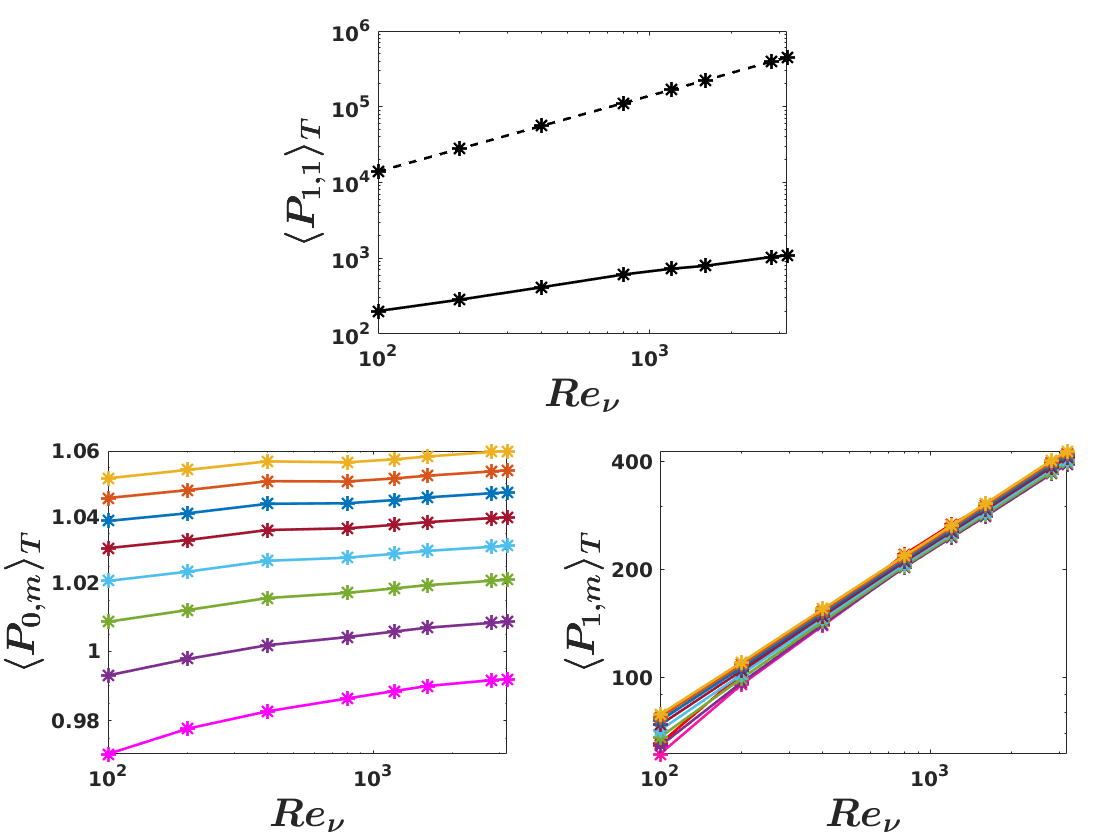}
		\put(-100,150){\large A}
		\put(90,150){\large B}
	
		\includegraphics[width=0.5\linewidth]{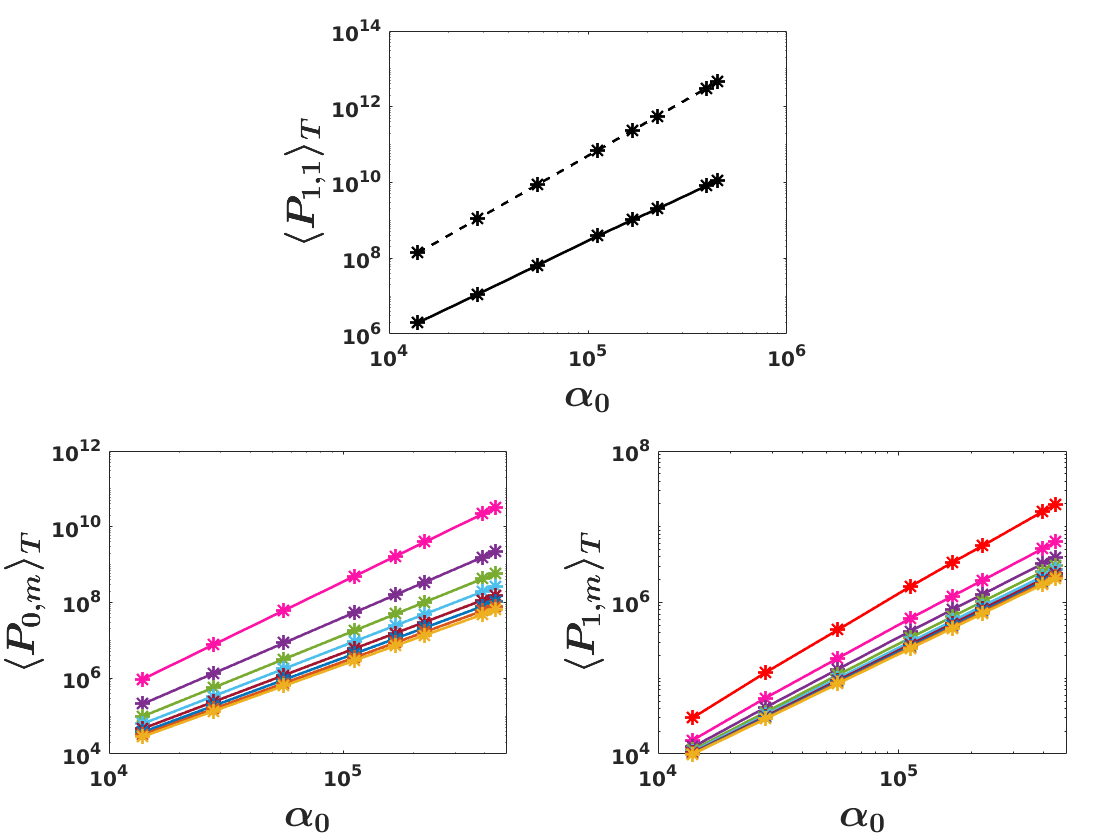}};
		\draw[line width=0.2mm,black, -] (0,0.0) -- (6.7,0.0);
		\draw[line width=0.2mm,black, -] (0,0.0) -- (0.0,5.2);
		\draw[line width=0.2mm,black, -] (6.7,0.0) -- (6.7,5.2);
		\draw[line width=0.2mm,black, -] (0,5.2) -- (6.7,5.2);
	
		\draw[line width=0.2mm,black, -] (6.75,0) -- (13.6,0);
		\draw[line width=0.2mm,black, -] (6.75,0) -- (6.75,5.2);
		\draw[line width=0.2mm,black, -] (6.75,5.2) --(13.6,5.2);
		\draw[line width=0.2mm,black, -] (13.6,0) --(13.6,5.2);
		\end{tikzpicture}

			\caption{(Colour online) Illustrative plots for $U_{0} = \sqrt{\alpha/\beta}$ (panel A) and $U_{0} = \nu/L$ (panel B) for $d=2$, runs A1-A8 (see Table \ref{tab:parameters})\,: First  row\,: plots versus $Re_{\nu}$ (panel A) and $\alpha_{0}$ (panel B) of $\left< P_{1,1} \right> _T$ ( solid black line) and $Re_{\nu}\,\alpha_{0}\,\mathcal{A}_{0}$ (dashed black line). Second row\,: Plots versus $Re_{\nu}$ (panel A) and $\alpha_{0}$ (panel B) of $\left<  P_{0,m}  \right> _T$ and $\left<P_{1,m}\right> _T$. Curves for
			$m=2, 3, 4, 5, 6, 7, 8, 9$, and $10$ are drawn in red, pink; violet, green, cyan, maroon, blue, orange, and yellow, respectively.
			Similar plots for other representative runs are given in the Supplemental Material.}
				\label{fig:fig2}
		\end{figure*}
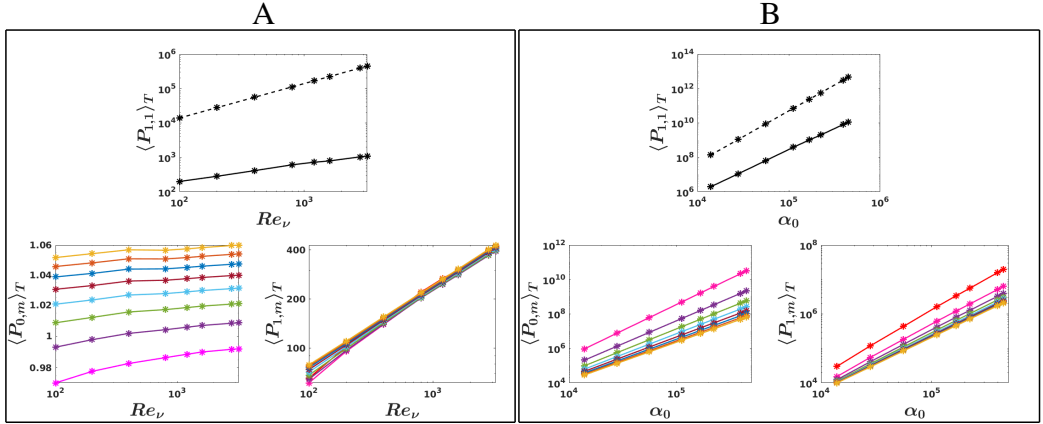

\smallskip
{In Fig.~\ref{fig:fig2} we present plots for} $U_{0} = \sqrt{\alpha/\beta}$ (panel A) and $U_{0} = \nu/L$ (panel B), runs A1-A8, to illustrate whether the bound in \ref{ee4a} is saturated\,: In the first row we plot $\left< P_{1,1} \right> _T$ (solid black line) versus $Re_{\nu}$ (panel A) and $\alpha_{0}$ (panel B)\,; the black dashed line denotes $Re_{\nu}\,\alpha_{0}\,\mathcal{A}_{0}$, which is the right-hand side (RHS) of (\ref{ee4a}) . In the second row we present plots versus $Re_{\nu}$ (panel A) and $\alpha_{0}$ (panel B) of $\left<P_{0,m}\right> _T$ and $\left<P_{1,m}\right> _T$, for $m=2, \, \ldots \,, 10$. Note that curves for $ \left<  P_{1,m}\right> _T$ can cross as $Re_{\nu}$ increases (panel A)\,; by contrast, they do not cross as $\alpha_{0}$ increases (panel B).  Similar plots for other representative runs are given in the Supplemental Material.

\par\medskip\noindent

\section{Summary of results in the $d=3$ case}\label{sec:sumd3}

The proof of the results in the following subsections are given in \ref{appB}. The methods used there are based on the differential inequalities explained in \ref{appinequal}.

\subsection{Estimates for $\left<Q_{n,m}\right>_{T}$}\label{subsec:est3}

Results in the $d=3$ case are more restricted, which reflects the open status of the regularity problem. Nevertheless, time averages of various $Q_{n,m}$ of Navier-Stokes type can be found \cite{JDG2019,JDG2020}. In addition to a bound on $\left<Q_{1,1}\right>_{T}$\,, as in (\ref{ee4b}), our results from \ref{appB} are summarised thus\,: from (\ref{f10}) we have
\bel{Q1a}
\left<Q_{2,1}\right>_{T} \leq c\, \alpha_{0}\Rn^{2}.
\ee 
We also find that for $n\geq 2$ and $m \geq 1$,
\bel{Q1b}
\left<Q_{n,m}\right>_{T} < \infty\,
\ee
although estimating the right hand side is a difficult calculation that we have omitted (see (\ref{Qn3})).
\par\smallskip\noindent
Moreover, with 
\bel{Q2a}
Q_{0,m} = \|\bu\|_{2m}^{\frac{2m}{2m-3}}\,,
\ee
{for $m > 2$, (\ref{um5}) shows that}
\bel{Q2b}
\left<Q_{0,m}\right>_{T} \leq c\,\Ra^{\frac{2(m+3)}{5(2m-3)}}
\left(\alpha_{0}\Rn^{2}\right)^{\frac{9(m-2)}{5(2m-3)}}\,.
\ee
(\ref{um6}) also shows that, in the limit $m\to\infty$, 
\beq{Q2c}
\left<\|\bu\|_{\infty}\right>_{T} &\leq&  c\,\alpha_{0}^{9/10}\Ra^{1/5}\Rn^{9/5}\,.
\eeq

\subsection{Numerical results for $Q_{n,m}$}\label{subsec:num3}

\begin{figure*}[!]

		
		\begin{tikzpicture}
		\node[anchor=south west,inner sep=0] at (0,0)
		{\includegraphics[width=0.5\linewidth]{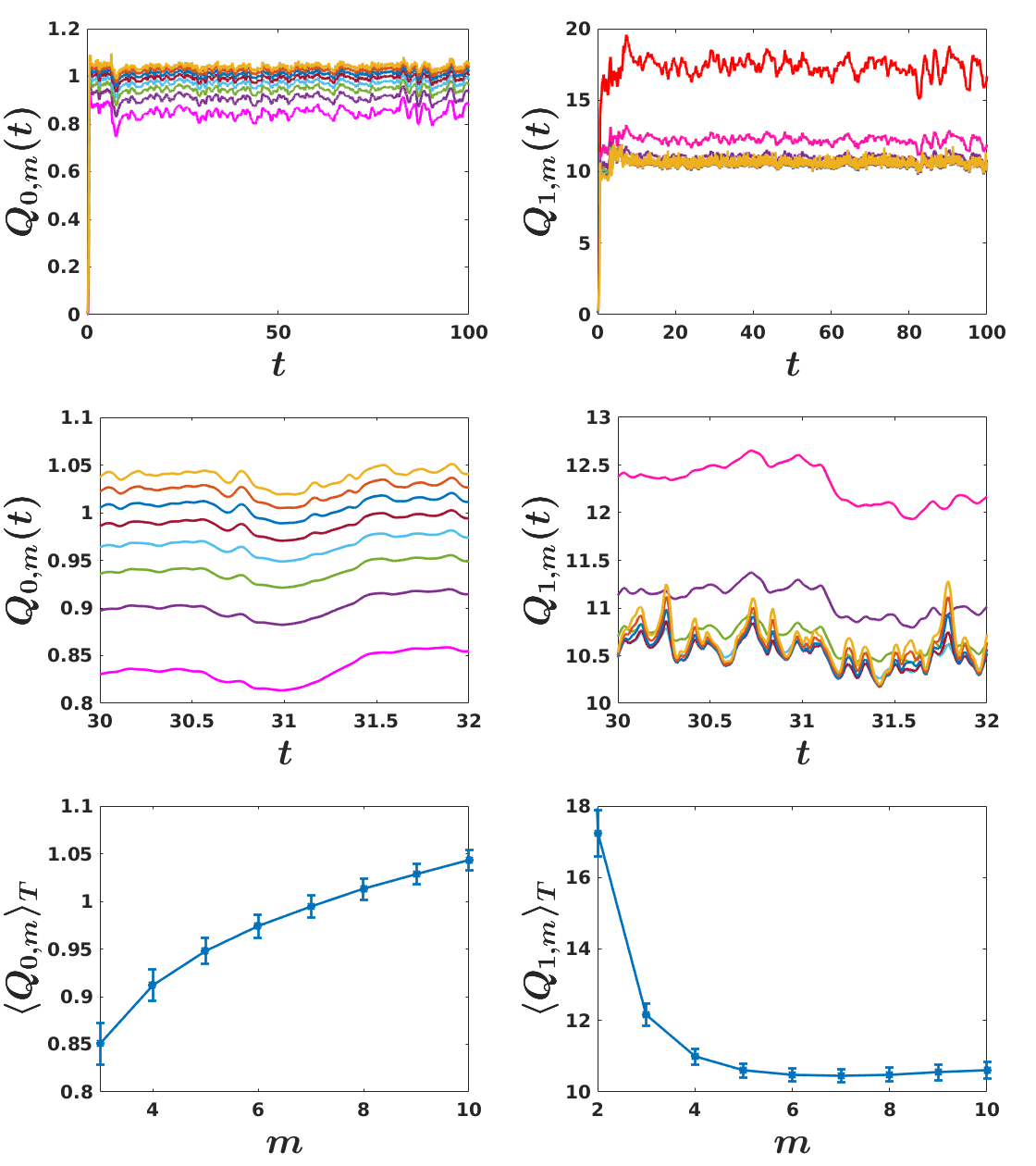}
		\put(-100,225){\large A}
		\put(90,225){\large B}
	
		\includegraphics[width=0.5\linewidth]{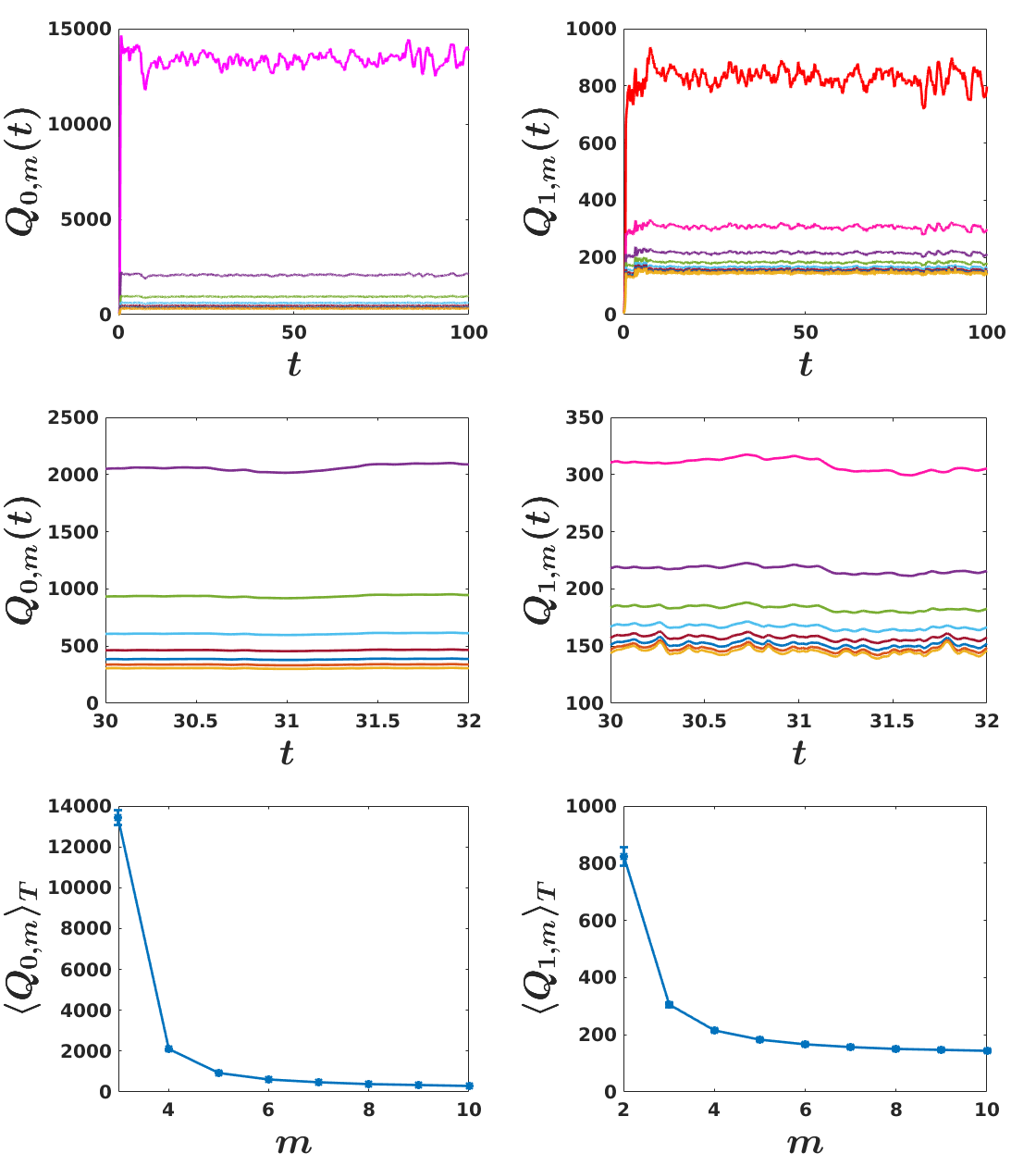}};
		\draw[line width=0.2mm,black, -] (0,0.0) -- (6.7,0.0);
		\draw[line width=0.2mm,black, -] (0,0.0) -- (0.0,7.8);
		\draw[line width=0.2mm,black, -] (6.7,0.0) -- (6.7,7.8);
		\draw[line width=0.2mm,black, -] (0,7.8) -- (6.7,7.8);
		\draw[line width=0.1mm,blue, dashed] (1.825,5.55) -- (0.7,4.9);
		\draw[line width=0.1mm,blue, dashed] (2.0,5.55) -- (3.1,4.9);
		\draw[line width=0.1mm,blue, dashed] (5.2,5.55) -- (4.0,4.9);
		\draw[line width=0.1mm,blue, dashed] (5.3,5.55) -- (6.4,4.9);
		\draw[line width=0.2mm,black, -] (6.75,0) -- (13.6,0);
		\draw[line width=0.2mm,black, -] (6.75,0) -- (6.75,7.8);
		\draw[line width=0.2mm,black, -] (6.75,7.8) --(13.6,7.8);
		\draw[line width=0.2mm,black, -] (13.6,0) --(13.6,7.8);
		\draw[line width=0.1mm,blue, dashed] (8.625,5.55) -- (7.5,4.9);
		\draw[line width=0.1mm,blue, dashed] (8.8,5.55) -- (9.9,4.9);
		\draw[line width=0.1mm,blue, dashed] (11.9,5.55) -- (10.8,4.9);
		\draw[line width=0.1mm,blue, dashed] (12.1,5.55) -- (13.2,4.9);
		\end{tikzpicture}
			\caption{(Colour online) Illustrative plots for $U_{0} = \sqrt{\alpha/\beta}$ (panel A) and $U_{0} = \nu/L$ (panel B) for $d=3$, run B2 (see Table \ref{tab:parameters}): First and second rows: plots versus $t$ of $Q_{0,m}$ and $Q_{1,m}$; the plots in the second row are expanded versions of small segments of the plots in the first row. Third row: Plots versus $m$ of $\left< Q_{0,m} \right>_T$ and $\left< Q_{1,m} \right>_T$. Curves for
			$m=2, 3, 4, 5, 6, 7, 8, 9$, and $10$ are drawn in red, pink; violet, green, cyan, maroon, blue, orange, and yellow, respectively.
			Similar plots for other representative runs are given in the Supplemental Material.}
			\label{fig:fig3}
				
		\end{figure*}
		\begin{figure*}[!]


		
		\begin{tikzpicture}
		\node[anchor=south west,inner sep=0] at (0,0)
		{\includegraphics[width=0.5\linewidth]{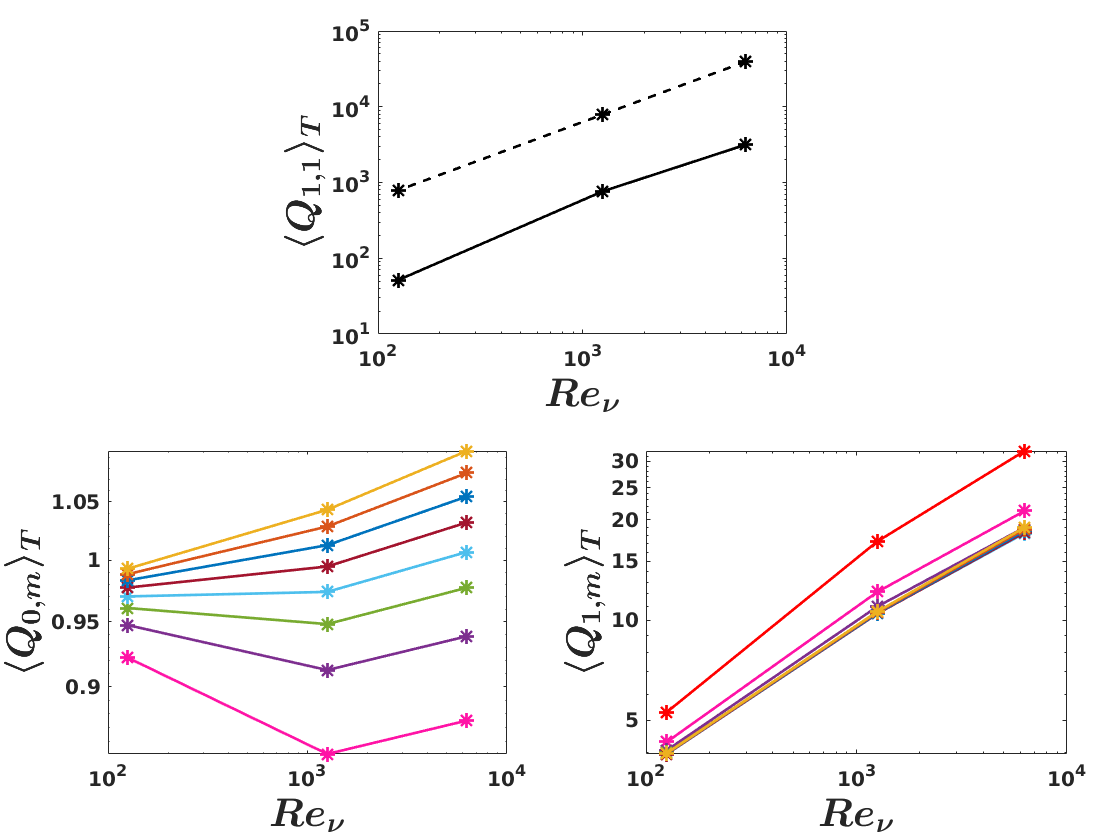}
		\put(-100,150){\large A}
		\put(90,150){\large B}
	
		\includegraphics[width=0.5\linewidth]{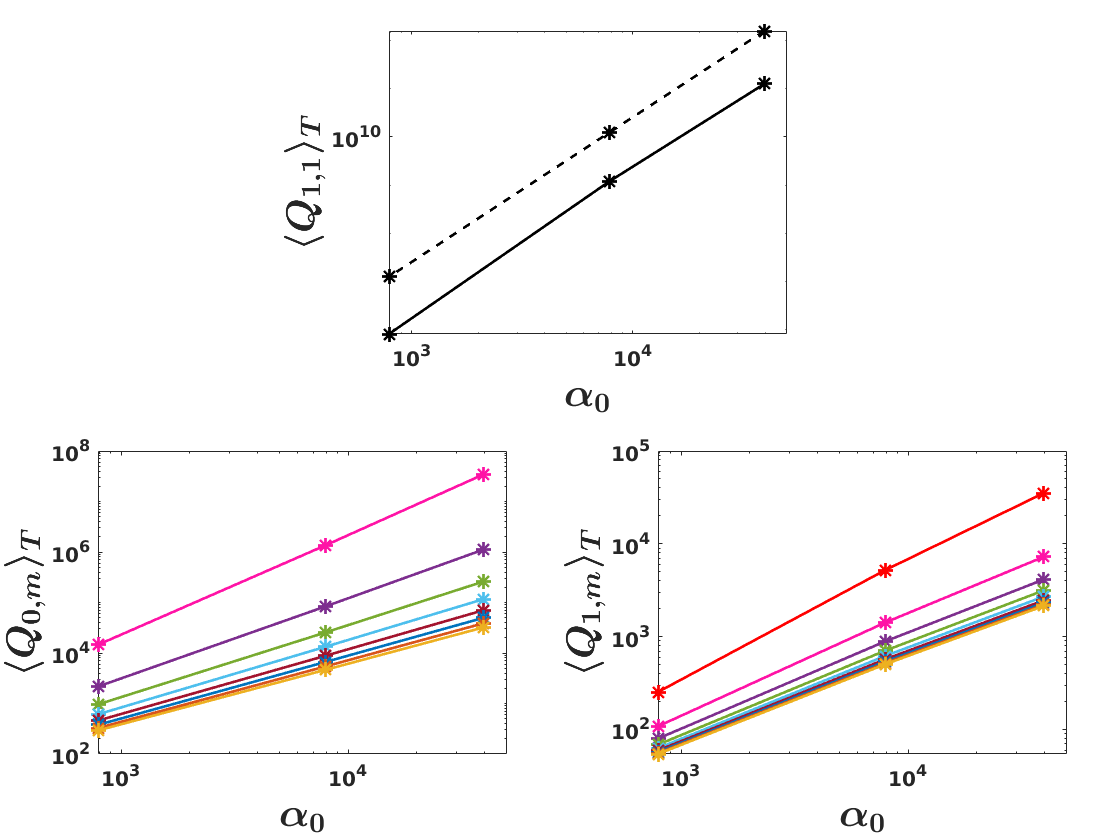}};
		\draw[line width=0.2mm,black, -] (0,0.0) -- (6.7,0.0);
		\draw[line width=0.2mm,black, -] (0,0.0) -- (0.0,5.2);
		\draw[line width=0.2mm,black, -] (6.7,0.0) -- (6.7,5.2);
		\draw[line width=0.2mm,black, -] (0,5.2) -- (6.7,5.2);
	
		\draw[line width=0.2mm,black, -] (6.75,0) -- (13.6,0);
		\draw[line width=0.2mm,black, -] (6.75,0) -- (6.75,5.2);
		\draw[line width=0.2mm,black, -] (6.75,5.2) --(13.6,5.2);
		\draw[line width=0.2mm,black, -] (13.6,0) --(13.6,5.2);
		
		\end{tikzpicture}
	\caption{(Colour online) Illustrative plots for $U_{0} = \sqrt{\alpha/\beta}$ (panel A) and $U_{0} = \nu/L$ (panel B) for $d=3$, runs A1-A8 (see Table \ref{tab:parameters}): First  row\,: plots versus $Re_{\nu}$ (panel A) and $\alpha_{0}$ (panel B) of $\left< Q_{1,1} \right> _T$  ( solid black line) and $Re_{\nu}\alpha_{0}\mathcal{A}_{0}$ (dashed black line). Second row\,: Plots versus $Re_{\nu}$ (panel A) and $\alpha_{0}$ (panel B) of $ \left<Q_{0,m}  \right> _T$ and $ \left<Q_{1,m}\right> _T$. Curves for
			$m=1, 2, 3, 4, 5, 6, 7, 8, 9$, and $10$ are drawn in black, red, pink, violet, green, cyan, maroon, blue, orange, and yellow, respectively.
			Similar plots for other representative runs are given in the Supplemental Material.}
			\label{fig:fig4}
			\end{figure*}
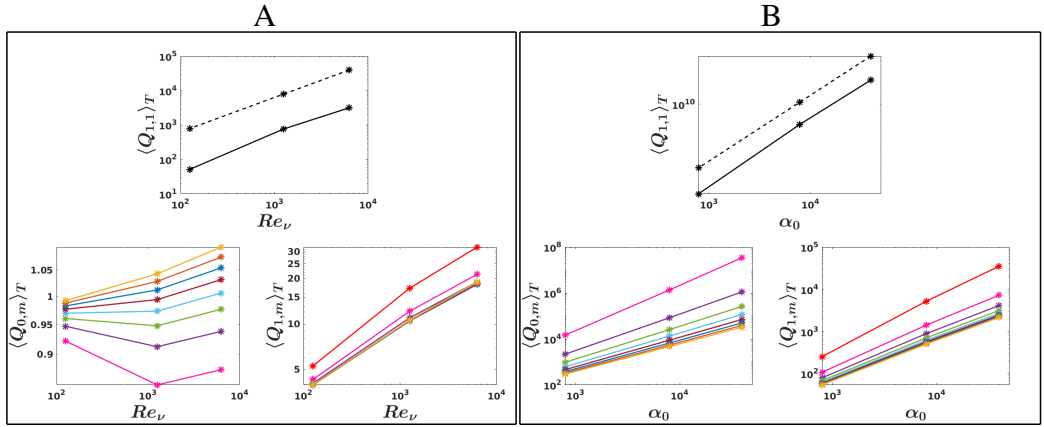

Fig. \ref{fig:fig3} shows plots of $Q_{n,m}(t)$ versus time $t$ for the two values $n=0$ and $n=1$, with a sequence of values of $m = 1,\, \ldots,\, 10$. Again, as in the $d=2$ case, there are two regimes, namely, an increasing regime $Q_{n,m} \leq Q_{n,m+1}$ and a decreasing regime $Q_{n,m}\geq Q_{n,m+1}$, because the norms $\|\cdot\|_{2m}$ must be ordered with increasing $m$, such that $\|\cdot\|_{2m} \leq \|\cdot\|_{2(m+1)}$; but the $\alpha_{n,m,d=3}$ \textit{decrease} as $m$ increases. 

{When $U_{0} = \sqrt{\alpha/\beta}$ (panel A)} we have observed both these regimes, but only the decreasing regime $Q_{n,m} \geq Q_{n,m+1}$ when $U_{0} = \nu/L$ (panel B) has been used. We illustrate this with plots in Fig.~\ref{fig:fig3} for run B2\,: In panel A ($U_{0} = \sqrt{\alpha/\beta}$) the plots of $Q_{1,m}$ can cross each other at different times, as we can see clearly in the expanded plots in the second row; such crossings do not occur in panel B ($U_{0} = \nu/L$). Furthermore, the plots versus $m$ of $\left< Q_{1,m} \right>_{T}$ (third row) decreases monotonically with increasing $m$ in panel B but not in panel A. Note, in particular, that $Q_{1,m}$ is almost equivalent to a non-dimensionalized version of the $D_{m}=\|\bom\|^{\alpha_{n,m,d=3}}_{2m}$ introduced in \cite{DKGGPV2013} and plotted there for the $d=3$ NSEs (see the Supplemental Material). The only difference here is that we are plotting $\|\nabla\bu\|_{2m}$ and not $\|\bom\|_{2m}$\,: the two are identical only when $m=1$.
\par\smallskip
In Fig.~\ref{fig:fig4} we present plots for $U_{0} = \sqrt{\alpha/\beta}$ (panel A) and $U_{0} = \nu/L$ (panel B), runs B1-B3, to illustrate whether the bound in (\ref{ee4a}) is saturated\,: In the first row we plot $\left< Q_{1,1} \right> _T$ (panel B)\,; the black dashed line denotes $Re_{\nu}\,\alpha_{0}\,\mathcal{A}_{0}$, which is the RHS of (\ref{ee4a}). In the second row we present plots versus $Re_{\nu}$ (panel A) and $\alpha_{0}$ (panel B) of $\left<  Q_{0,m}\right> _T$ and $ \left<Q_{1,m}\right> _T$, for $m=2, \, \ldots \,, 10$. Note that curves for $ \left<Q_{1,m}  \right> _T$ can cross as $Re_{\nu}$ increases (panel A)\,; by contrast, they do not cross as $\alpha_{0}$ increases (panel B). Similar plots for other representative runs are given in the Supplemental Material.

\section{Conclusions}\label{sec:Conclusions}

We have married the two approaches of the analysis of solutions of the ITT equations through the estimation of weighted time averages, together with the results of numerical simulations. To achieve this we have invoked the similar scaling properties between the ITT equations and the NSEs\,: see \S\ref{sec:scaling}. There are, however, two important differences. Usually the NSEs are considered either in the decaying or the additively forced case, whereas the ITT equations have no additive forcing but instead possess a linear-activity term $\alpha_{0}\bu$ which, in effect, pumps energy into the system. {Dynamically the effect of this term, together with the negative cubic term, creates a platform for either temporally frozen solutions or statistically steady states. Furthermore, it is shown in \S\ref{subsec:reg1} that the greatest contrast with the $2d$ NSEs lies at the level of the vorticity equation. Whereas the absence of vortex stretching in $2d$ NSEs lies at the root of the regularity of its solutions, the $\mbox{curl}\left(\bu|\bu|^{2}\right)$-term in the ITT equations appears to recreate another form of vortex stretching. However, the extra piece of information afforded to us is the bounded time average of the $L^{4}$ norm of the velocity field expressed in (\ref{ee3c}), which is just enough to recover regularity, but only to the degree that bounds are exponential in time -- see (\ref{reg1}). Thus, we fall just short of proving the existence of a \textit{global attractor}. Results that lie in parallel with those of the NSEs in both spatial dimensions are the existence of bounded infinite hierarchies of time averages\,: i.e., estimates of   $\left<P_{n,m}\right>_{T}$ and $\left<Q_{n,m}\right>_{T}$\,, whose bounds are calculated in \ref{appA} and \ref{appB} and summarised in \S\ref{sec:sumd2} and \S\ref{sec:sumd3}, together with the results of our numerical simulations. When statistically steady solutions appear, the possibility of multifractality should be considered \cite{RP2020,RP2022}. A future line of approach might be to repeat the calculation in \cite{BD-JDG}, where the correspondence between the multifractal model of turbulence and the NSEs was investigated.}
\par\bigskip\noindent
\textit{Acknowledgments\,:} RP and JDG would like to thank the Isaac Newton Institute for Mathematical Sciences for support and hospitality  during the programme \textit{Mathematical Aspects of Fluid Turbulence\,: where do we stand?} in 2022, when work on this paper was undertaken. It was supported by grant number EP/R014604/1. RP, KVK and NBP thank SERB (India), the National Supercomputing Mission (India), and SERC (IISc), for support and computational resources, and A. Gupta for discussions.

\appendix
\section{Differential Inequalities}\label{appinequal}

The most widely used class of differential inequalities that generalize the Sobolev inequalities are called Gagliardo-Nirenberg inequalities \cite{Adams1975}. In their most general form in integer $d$ dimensions ($d=1,2,3$) they can be expressed as 
\bel{GN1}
\|\nabla^{j}u\|_{p} \leq c\, \|\nabla^{n}u\|_{r}^{a}\|u\|_{q}^{1-a}\,,
\ee
where $0 \leq j < n$ and $1 < p,\,r,\,q \leq \infty$. The exponent $a$ can be calculated by dimensional analysis and thus must satisfy 
\bel{GN2}
\frac{1}{p} =  \frac{j}{d} + a\left(\frac{1}{r} - \frac{n}{d}\right) + \frac{1-a}{q}\,,
\ee
{where $j/n \leq a < 1$. (\ref{GN2}) holds on the whole space $\mathbb{R}^{d}$. With periodic boundary conditions, there are $L^{2}$ additional terms there  are lower-order corrections to our estimates, which will be ignored. In the following\,:}
\ben\itemsep 0mm
\item {The standard Leray-Hopf weak-solution machinery has been used, which has been derived for the NSEs in \cite{FGT1981}.} 

\item We will use the convention that the constants designated as $c$, $c_{m}$ or $c_{n,m}$ are generic in the sense that they may differ from line to line. 
\een

\section{Proofs in the $d=2$ case}\label{appA}

\subsection{Estimates for $\left<P_{2,1}\right>_{T}$}\label{P2sect}

Our first requirement is to bound $P_{2,1}$, with $P_{n,1}$ defined in (\ref{Pnmdef}). Clearly $\alpha_{n,1,2} = 1/n$, so
\beq{P2}
\left<P_{2,1}\right>_{T} = \left<H_{2}^{1/2}\right>_{T} &=& \left<\left(\frac{H_{2}}{H_{1}}\right)^{1/2}H_{1}^{1/2}\right>_{T}\nonumber\\
&\leq& \left<\frac{H_{2}}{H_{1}}\right>_{T}^{1/2}
\left<P_{1,1}\right>_{T}^{1/2}\,.
\eeq
To achieve a bound on $\left<H_{2}/H_{1}\right>_{T}$\,, we take the curl of the ITT equation\,:
\bel{P3}
\left(\partial_{t} + \bu\cdot\nabla\right)\bom = \alpha_{0}\bom + \Rn^{-1}\Delta\bom - \Rb\,\mbox{curl}(\bu|\bu|^{2})\,.
\ee
The key point is that, while the vortex-stretching term $\bom\cdot\nabla\bu$ is missing because of the orthogonality of $\bom$ with the gradient operator, there is an additional negative $\Rb\,\mbox{curl}(\bu|\bu|^{2})$-term. To deal with this we write {
\beq{P4}
\shalf \dot{H}_{1} &=& \alpha_{0}H_{1} - \Rn^{-1}H_{2} - \Rb\I \bom\cdot \mbox{curl}(\bu|\bu|^{2})dV\nonumber\\
&\leq& \alpha_{0}H_{1} - \shalf \Rn^{-1}H_{2} + \shalf \Rb^{2}\Rn\I |\bu|^{6}dV\,,
\eeq
where we have integrated by parts and have then used a H\"older inequality.} Now divide by $H_{1}$ to obtain 
\bel{P5}
\shalf \left<\frac{H_{2}}{H_{1}}\right>_{T} \leq \Rn\alpha_{0} + \shalf\Rb^{2}\Rn^{2}\left<\frac{\I |\bu|^{6}dV}{H_{1}}\right>_{T}\,.
\ee
As in (\ref{reg3}), we use the Gagliardo-Nirenberg inequality
\bel{P6}
\|\bu\|_{6} \leq c\, \|\nabla\bu\|_{2}^{a}\|\bu\|_{4}^{1-a}\,, \qquad a = \sthird \,,
\ee
to find
\bel{P7}
\left<\frac{\|\bu\|_{6}^{6}}{H_{1}}\right>_{T} \leq c\,\left<\|\bu\|_{4}^{4}\right>_{T}\,.
\ee
{Inserting this into (\ref{P5}) and using (\ref{P6}) and the definition of $\Ra$ in (\ref{ee3c}), we find that
\beq{P8}
\shalf\left<\frac{H_{2}}{H_{1}}\right>_{T} &\leq& \alpha_{0}\Rn + c\,\Rb^{2}\Rn^{2}\left<\|\bu\|_{4}^{4}\right>_{T}\nonumber\\
&\leq& c\,\alpha_{0}\Rn\left(1 + \alpha_{0}\Rn\right)\,.
\eeq
Thus, to leading order in $\Rn$, (\ref{P2}) becomes 
\bel{P9}
\left<P_{2,1}\right>_{T} \lessapprox c\,\Ra^{1/2}\left(\alpha_{0}\Rn\right)^{3/2}\,,
\ee
as advertised in (\ref{P1c}). 
\par\smallskip
{For the second option in which $U_{0}$ is chosen to be $U_{0}= \nu/L$, thus making $\Rn=1$, a re-working of the bounds above gives 
\bel{P9extra}
\left<P_{2,1}\right>_{T} \leq c\, \alpha_{0}^{3/2}\left(1+\alpha_{0}\right)^{1/2}\Rb^{-1/2}\,.
\ee
We leave the remaining members of this second class of estimates to be calculated by the reader.}

\subsection{An estimate for $P_{1,m} =\left<\|\nabla\bu\|_{2m}^{\frac{2m}{2m-1}}\right>_{T}$}
\label{P1m}

A Gagliardo-Nirenberg inequality shows that, for some function $A$, for $m \geq 1$
\bel{P1ma}
\|A\|_{2m} \leq c_{m}\|\nabla A\|_{2}^{\frac{m-1}{m}}\|A\|_{2}^{\frac{1}{m}}\,.
\ee
We choose $A=\nabla\bu$ and, noting from (\ref{dim2a}) that $\alpha_{1,m,2}= m/(2m-1)$, we write
\beq{P1mb}
\left<\|\nabla\bu\|_{2m}^{\frac{2m}{2m-1}}\right>_{T} &\leq& c_{m}
\left<\|\nabla^{2}\bu\|_{2}^{\frac{2(m-1)}{2m-1}}\|\nabla\bu\|_{2}^{\frac{2}{2m-1}}\right>_{T}\nonumber\\
&\leq& c_{m}\left<\|\nabla^{2}\bu\|_{2}\right>_{T}^{\frac{2(m-1)}{2m-1}}\left<\|\nabla\bu\|_{2}^{2}\right>_{T}^{\frac{1}{2m-1}}\nonumber\\
&=& c_{m}\left<P_{2,1}\right>_{T}^{\frac{2(m-1)}{2m-1}}\left<H_{1}\right>_{T}^{\frac{1}{2m-1}}\nonumber\\
&\leq& c_{m}\Ra^{\frac{m}{2m-1}}\left(\alpha_{0}\Rn\right)^{\frac{3m-2}{2m-1}}\,,
\eeq
as advertised in (\ref{P2c}). In the limit $m\to\infty$\,,
\bel{P1mc}
\left<\|\nabla\bu\|_{\infty}\right>_{T} \leq c\,\Ra^{1/2}\left(\alpha_{0}\Rn\right)^{3/2}\,.
\ee
}
\subsection{Estimates for $\left<P_{0,m}\right>_{T} = \left<\|\bu\|_{2m}^{\frac{2m}{m-1}}\right>_{T}$ and $\left<\|\bu\|_{\infty}^{2}\right>_{T}$}\label{uinf}

We now turn to estimating $\bu$ in $L^{2m}(V)$ for $m > 2$. A Gagliardo-Nirenberg inequality shows that 
\bel{P10}
\|\bu\|_{2m}\leq c_{m}\|\nabla^{2}\bu\|_{2}^{a}\|\bu\|_{4}^{1-a}\,,
\ee
where $a = \frac{m-2}{3m}$. When $n=0$ and $d=2$ we have $(4-d)\alpha_{0,m,2} = \frac{2m}{m-1}$. Thus, 
\beq{P11}
\left<\|\bu\|_{2m}^{\frac{2m}{m-1}}\right>_{T} &\leq& c_{m}\left<\|\nabla^{2}\bu\|_{2}^{\frac{2(m-2)}{3(m-1)}}
\|\bu\|_{4}^{\frac{4(m+1)}{3(m-1)}}\right>_{T}\nonumber\\
&\leq& \left<P_{2,1}\right>_{T}^{\frac{2(m-2)}{3(m-1)}}\left<\|\bu\|_{4}^{4}\right>_{T}^{\frac{(m+1)}{3(m-1)}}\,,
\eeq
in which case, for $m > 2$, using (\ref{P9}) and (\ref{ee3c}),
\bel{P12}
\left<\|\bu\|_{2m}^{\frac{2m}{m-1}}\right>_{T} 
\leq c\,\Ra^{\frac{m}{m-1}}\left(\alpha_{0}\Rn\right)^{\frac{m-2}{m-1}}\,,
\ee
as advertised in (\ref{P2e}). In the limit $m\to\infty$\,, this reduces to 
\bel{P13}
\left<\|\bu\|_{\infty}^{2}\right>_{T} \leq c\,\alpha_{0}\Ra\Rn\,,
\ee
as advertised in (\ref{P2f}). 


\subsection{Estimates for $\left<P_{n,1}\right>_{T}$ for $n > 2$}\label{Pnsect}

Using the methods in \cite{DG1995}, a full `ladder' for $H_{n}$ takes the form
\beq{pn1}
\shalf \dot{H}_{n} &\leq& \alpha_{0}H_{n} - \Rn^{-1}H_{n+1} + c_{n,1}H_{n+1}^{1/2}H_{n}^{1/2}\|\bu\|_{\infty}\nonumber\\
&+& c_{n,2}\Rb H_{n}\|\bu\|_{\infty}^{2}\,;
\eeq
therefore, after a H\"older inequality and re-arrangement, we have 
\bel{pn2}
\shalf \dot{H}_{n}\leq \alpha_{0}H_{n} - \shalf\Rn^{-1}H_{n+1} + c_{n}\left(\Rb + \Rn\right)H_{n}\|\bu\|_{\infty}^{2}\,.
\ee
Thus,
\bel{pn3}
\left<\frac{H_{n+1}}{H_{n}}\right>_{T} \leq 2\alpha_{0}\Rn + c_{n}\Rn\left(\Rb + \Rn\right)\left<\|\bu\|_{\infty}^{2}\right>_{T}\,.
\ee
{From this and (\ref{P13}) we deduce that }
\beq{pn4}
\left<\frac{H_{n+1}}{H_{n}}\right>_{T} 
&\lessapprox& c_{n}\alpha_{0}\Ra\Rn^{3}\,.
\eeq
Moreover, for $n\geq 2$
\beq{pn5}
\left<P_{n+1,1}\right>_{T} &=& \left<\left(\frac{H_{n+1}}{H_{n}}\right)^{\frac{1}{n+1}}H_{n}^{\frac{1}{n+1}}\right>_{T}\nonumber\\
&\leq& \left<\frac{H_{n+1}}{H_{n}}\right>_{T}^{\frac{1}{n+1}}\left<P_{n,1}\right>_{T}^{\frac{n}{n+1}}\,.
\eeq
{Given that $\left<P_{2,1}\right>_{T}$ is bounded above as in (\ref{P9}), then together with (\ref{pn4}), we can now generate upper bounds on every $\left<P_{n,1}\right>_{T}$. For $n \geq 2$, these are}
\beq{Pn6}
\left<P_{n,1}\right>_{T} &\leq& c_{n,1}\alpha_{0}^{\frac{n+1}{n}}\Ra^{\frac{n-1}{n}}\Rn^{\frac{3(n-1)}{n}}\,,
\eeq
which can be transformed into the form advertised in (\ref{P2a}).


\subsection{Estimates for $\left<P_{n,m}\right>_{T}$}\label{Pnmsect}

We can write down an inequality for $B = \nabla^{2}\bu$ such that
\bel{Pnms1}
\|\nabla^{n-2}B\|_{2m} \leq c\,\|\nabla^{N-2}B\|_{2}^{a}\|B\|_{2}^{1-a}\,,
\ee
for some $N > n + 1 - 1/m$\,, with
\bel{Pnms2}
a = \frac{m(n-1)-1}{m(n+1)-1}\,.
\ee
Thus, we can write 
\bel{Pnms3}
\left<\|\nabla^{n}\bu\|_{2m}^{2\alpha_{n,m}}\right>_{T} 
\leq c\,\left<\|\nabla^{N}\bu\|_{2}^{2a\alpha_{n,m}}\|\nabla^{2}\bu\|_{2}^{2(1-a)\alpha_{n,m}}\right>_{T}\,.
\ee
Re-arranging and then using H\"older's inequality, we have
\beq{Pnms4}
\left<\|\nabla^{n}\bu\|_{2m}^{2\alpha_{n,m}}\right>_{T} &\leq& 
c\,\left<\left(\|\nabla^{N}\bu\|_{2}^{2/N}\right)^{aN\alpha_{n,m}}\left(\|\nabla^{2}\bu\|_{2}\right)^{2(1-a)\alpha_{n,m}}\right>_{T}\\
&\leq& 
c_{N,n,m}\left<\|\nabla^{N}\bu\|_{2}^{2/N}\right>_{T}^{aN\alpha_{n,m}}\left<\|\nabla^{2}\bu\|_{2}^{\frac{2(1-a)\alpha_{n,m}}{1-aN\alpha_{n,m}}}\right>_{T}^{1-aN\alpha_{n,m}}\\
&=& 
c_{N,n,m}\left<P_{N,1}\right>_{T}^{aN\alpha_{n,m}}\left<P_{2,1}^{\frac{2(1-a)\alpha_{n,m}}{1-aN\alpha_{n,m}}}\right>_{T}^{1-aN\alpha_{n,m}}\,.
\eeq
It can be checked through the definition of $a$ in (\ref{Pnms2}) that the exponent of $P_{2,1}$ inside the time-average is unity. Estimates for $\left<P_{N,1}\right>_{T}$ and $\left<P_{2,1}\right>_{T}$ are available from (\ref{P2a}) and (\ref{P1c})\,: one can then choose the lowest value of $N$\,, constrained by $N > n + 1 - 1/m$. After some algebra, this leads to the result
\bel{Pnms5}
\left<P_{n,m}\right>_{T} \leq c_{n,m}\alpha_{0}^{\frac{2m}{m(n+1)-1}}
\left(\alpha_{0}\Ra\Rn^{3}\right)^{\frac{mn-1}{m(n+1)-1}}\,,
\ee
as advertised in (\ref{P2g}).


\section{Proofs in the $d=3$ case}\label{appB}

\par\smallskip\noindent
\textbf{Step 1\,:}  $Q_{n,m}$ is defined in (\ref{Qnmdef}). In simplied form $Q_{n,1}$ can be written as
\bel{f1a}
Q_{n,1} = H_{n}^{\frac{1}{2n-1}}\,.
\ee
Moreover, because $\left<Q_{1,1}\right>_{T} = \left<H_{1}\right>_{T}$\,, we have an estimate for this in (\ref{ee3b}). We begin this section by estimating $\left<Q_{2,1}\right>_{T}$ from the vorticity equation (\ref{P3}), with the $3d$ vortex stretching term restored\,:
\bel{f1b}
\left(\partial_{t} + \bu\cdot\nabla\right)\bom = \alpha_{0}\bom + \Rn^{-1}\Delta\bom + \bom\cdot\nabla\bu - \Rb\,\mbox{curl}(\bu|\bu|^{2})\,.
\ee
{The equivalent of (\ref{P4}) is
\beq{f2}
\shalf \dot{H}_{1} &\leq& \alpha_{0}H_{1} - \Rn^{-1}H_{2} + c\,H_{2}^{1/2}H_{1}^{1/2}\|\bu\|_{\infty} 
- \Rb\I \bom\cdot \mbox{curl}(\bu|\bu|^{2})dV\nonumber\\
&\leq& \alpha_{0}H_{1} - \sthreequart \Rn^{-1}H_{2} + c\,H_{2}^{1/2}H_{1}^{1/2}\|\bu\|_{\infty}  + \Rb^{2}\Rn\I |\bu|^{6}dV\,,
\eeq
where we have integrated by parts and have then used a H\"older inequality.} The $3$-dimensional Agmon inequality for $\|\bu\|_{\infty}$ is 
\bel{f3}
\|\bu\|_{\infty}^{2}\leq c_{n}H_{n}^{a}H_{1}^{1-a}\qquad n \geq 2\,,
\ee
with $a=\frac{1}{2(n-1)}$. Thus, for $n=2$\,,
\bel{f4}
\|\bu\|_{\infty} \leq c\,H_{2}^{1/4}H_{1}^{1/4}
\ee
and so 
\beq{f5}
H_{2}^{1/2}H_{1}^{1/2}\|\bu\|_{\infty} &\leq& c\,H_{2}^{3/4}H_{1}^{3/4}\\
&\leq& \squart \Rn^{-1}H_{2} + {\sthreequart}c\,\Rn^{3} H_{1}^{3}\,.
\eeq
Moreover, Sobolev's inequality for $d=3$ shows that 
\bel{f6}
\|\bu\|_{6} \leq c\,\|\nabla\bu\|_{2}\,.
\ee
{Therefore,} in total, (\ref{f2}) becomes 
\bel{f7a}
\shalf \dot{H}_{1} \leq \alpha_{0}H_{1} - \shalf \Rn^{-1}H_{2} + c\,\Rn\left(\Rb^{2} + \Rn^{2}\right)H_{1}^{3}\,.
\ee
Thus, the ultimate contribution to (\ref{f7b}) from the $\bu|\bu|^{2}$-term is proportional to that from the vortex-stretching term, in the sense that they are both proportional to $H_{1}^{3}$.
Dividing by $H_{1}^{2}$ gives
\bel{f7b}
{\frac{1}{2}}\Rn^{-1}\left<\frac{H_{2}}{H_{1}^{2}}\right>_{T} 
\leq \alpha_{0} \left<H_{1}^{-1}\right>_{T} + c\,\Rn\left(\Rb^{2} + \Rn^{2}\right)
\left<Q_{1,1}\right>_{T}\,.
\ee
Ignoring the first term on the right hand side with the negative exponent, we can write 
\bel{f8}
\left<\frac{H_{2}}{H_{1}^{2}}\right>_{T} \leq c\,\alpha_{0}\Rn^{2}\left(\Rb^{2} + \Rn^{2}\right)\,.
\ee
{Thus, we finally have}
\beq{f9}
\left<Q_{2,1}\right>_{T} &=& \left<\left(\frac{H_{2}}{H_{1}^{2}}\right)^{1/3}H_{1}^{2/3}\right>_{T}\\
&\leq& \left<\frac{H_{2}}{H_{1}^{2}}\right>_{T}^{1/3}\left<H_{1}\right>_{T}^{2/3}\\
&\leq& c\,\left(\alpha_{0}\Rn^{2}\left(\Rb^{2} + \Rn^{2}\right)\right)^{1/3}(\alpha_{0}\Rn)^{2/3}\\
&=& c\,\alpha_{0}\Rn^{4/3}\left(\Rb^{2} + \Rn^{2}\right)^{1/3}\,.
\eeq
{If $\Rn$ is dominant, the bound scales like} 
\bel{f10}
\left<Q_{2,1}\right>_{T} \lessapprox \alpha_{0}\Rn^{2} + O\left(\Rn^{4/3}\right)\,.
\ee
\par\medskip\noindent
\textbf{Step 2\,:}  Let us repeat (\ref{pn2}) by writing
\bel{fgt1a}
\shalf \dot{H}_{n}\leq \alpha_{0}H_{n} - \shalf\Rn^{-1}H_{n+1} + c_{n}\left(\Rb + \Rn\right)H_{n}\|\bu\|_{\infty}^{2}\,.
\ee
After re-arrangement and the use of Agmon's inequality, (\ref{fgt1a}) becomes 
\beq{fgt2}
\shalf \dot{H}_{n}\leq \alpha_{0}H_{n} - \shalf\Rn^{-1}H_{n+1} + c_{n}\left(\Rb + \Rn\right)H_{n}^{1+a}H_{1}^{1-a}\,.
\eeq
Dividing by $H_{n}^{2n/(2n-1)}$ and time averaging gives
\beq{fgt3a}
\shalf\Rn^{-1}\left<\frac{H_{n+1}}{H_{n}^{\frac{2n}{2n-1}}}
\right>_{T} &\leq& \alpha_{0}\left<H_{n}^{1-\frac{2n}{2n-1}}\right>_{T} + c_{n}\left(\Rb + \Rn\right)\left<H_{n}^{\frac{2n-1}{2(n-1)} - \frac{2n}{2n-1}}H_{1}^{\frac{2n-3}{2(n-1)}}\right>_{T}
\nonumber\\
&\leq& \alpha_{0}\left<H_{n}^{-\frac{1}{2n-1}}\right>_{T} + c_{n}\left(\Rb + \Rn\right)
\left<Q_{n}\right>_{T}^{\frac{1}{2(n-1)}}\left<Q_{1}\right>^{\frac{2n-3}{2(n-1)}}\,.
\eeq
The next step is to ignore the first term\footnote{The term with the negative exponent of $H_{n}$ on the RHS of (\ref{fgt3a}) is only out of control if $H_{n}$ temporarily becomes very small. In principle, this could be dealt with by adding a constant term to $H_{n}$ to provide the platform of a lower bound. We omit the details.}. Given that $\Rn$ is the dominant term, we write (\ref{fgt3a}) in the simplified form
\bel{fgt3b}
\left<\frac{H_{n+1}}{H_{n}^{\frac{2n}{2n-1}}}\right>_{T} \leq c_{n}\Rn^{2}
\left<Q_{n,1}\right>_{T}^{\frac{1}{2(n-1)}}\left<Q_{1,1}\right>_{T}^{\frac{2n-3}{2(n-1)}}\,.
\ee
We then study
\beq{Qn1}
\left<Q_{n+1,1}\right>_{T} = \left< H_{n+1}^{\frac{1}{2n+1}}\right>_{T} &=& \left<\left(\frac{H_{n+1}}{H_{n}^{\frac{2n}{2n-1}}}
\right)^{\frac{1}{2n+1}}H_{n}^{\frac{2n}{(2n+1)(2n-1)}}\right>_{T}\nonumber\\
&\leq& \left<\frac{H_{n+1}}{H_{n}^{\frac{2n}{2n-1}}}\right>_{T}^{\frac{1}{2n+1}}\left<Q_{n,1}\right>_{T}^{\frac{2n}{2n+1}}\,,
\eeq
in which case
\beq{Qn2}
\left<Q_{n+1,1}\right>_{T} &\leq& c_{n}\Rn^{\frac{2}{2n+1}}
\left<Q_{n,1}\right>_{T}^{\frac{2n}{2n+1} + \frac{1}{2(2n+1)(n-1)}}\left<Q_{1,1}\right>_{T}^{\frac{2n-3}{2(n-1)(2n+1)}}\nonumber\\
&=& c_{n}\Rn^{\frac{2}{2n+1}}\left<Q_{n,1}\right>_{T}^{\frac{(2n-1)^{2}}{2(2n+1)(n-1)}}
\left<Q_{1,1}\right>_{T}^{\frac{2n-3}{2(n-1)(2n+1)}}\,.
\eeq
Given that we have estimates for both $\left<Q_{1,1}\right>_{T}$ and $\left<Q_{2,1}\right>_{T}$\,, we can generate estimates for all $\left<Q_{n,1}\right>_{T}$ for $n\geq 3$. {Thus, we deduce that
\bel{Qn3}
\left<Q_{n,1}\right>_{T} < \infty\qquad n\geq 3\,.
\ee
The bound is messy so we simply register that the right hand side is finite.  Finally, the method used in \S\ref{Pnmsect} can be used to show that $\left<Q_{n,m}\right>_{T} < \infty$ for $m \geq 1$.}
\par\smallskip\noindent
\textbf{Step 3\,:} Now let us consider
\bel{um1}  
\|\bu\|_{2m} \leq c\,\|\nabla^{2}\bu\|_{2}^{a}\|\bu\|_{4}^{1-a}\,,
\ee
where $a= 3(m-2)/5m$ with $m > 2$. {Because $\alpha_{2,1,3} = \twothirds$ we can write 
\beq{um2}
\left<\|\bu\|_{2m}^{\alpha_{0,m}}\right>_{T} &\leq& c\,\left<\|\nabla^{2}\bu\|_{2}^{a\alpha_{0,m}}\|\bu\|_{4}^{(1-a)\alpha_{0,m}}\right>_{T}\nonumber\\
&=& c\,\left<Q_{2,1}^{3a\alpha_{0,m}/2}\left(\|\bu\|_{4}^{4}\right)^{\squart(1-a)\alpha_{0,m}}\right>_{T}\,,
\eeq
where $\alpha_{0,m} = \frac{2m}{2m-3}$. Then 
\bel{um3}
\left<\|\bu\|_{2m}^{\alpha_{0,m}}\right>_{T}\leq c\,\left<Q_{2,1}\right>_{T}^{3a\alpha_{0,m}/2}
\left<\left(\|\bu\|_{4}^{4}\right)^{\frac{\squart(1-a)\alpha_{0,m}}{1-3a\alpha_{0,m}/2}}
\right>_{T}^{1-3a\alpha_{0,m}/2}\,.
\ee
Given $a$ and $\alpha_{0,m}$, it can easily be checked that the exponent of $\|\bu\|_{4}^{4}$ inside the average is unity. Thus, because $\left<\|\bu\|_{4}^{4}\right>_{T} \leq c\,\Ra^{2}$ and with the help of (\ref{f10}), (\ref{um3}) becomes}
\beq{um4}
\left<\|\bu\|_{2m}^{\alpha_{0,m}}\right>_{T} &\leq& c\,(\alpha_{0}\Rn^{2})^{3a\alpha_{0,m}/2}
\Ra^{2-3a\alpha_{0,m}}\nonumber\\
&=& \alpha_{0}^{3a\alpha_{0,m}/2}\Ra^{2-3a\alpha_{0,m}}\Rn^{3a\alpha_{0,m}}\,.
\eeq 
In fact $3a\alpha_{0,m}/2 = \frac{9(m-2)}{5(2m-3)}$ and so $1- 3a\alpha_{0,m}/2 = \frac{m+3}{5(2m-3)}$, 
{whence
\bel{um5}
\left<\|\bu\|_{2m}^{\alpha_{0,m}}\right>_{T} \leq 
c\,\Ra^{\frac{2(m+3)}{5(2m-3)}}\left(\alpha_{0}\Rn^{2}\right)^{\frac{9(m-2)}{5(2m-3)}}\qquad m > 2\,,
\ee
as advertised in (\ref{Q2b}). In the limit $m\to\infty$\,, we find that 
\beq{um6}
\left<\|\bu\|_{\infty}\right>_{T} &\leq& c\,\Ra^{1/5}\left(\alpha_{0}\Rn^{2}\right)^{9/10}\,,
\eeq
as advertised in (\ref{Q2c}).}



\bibliographystyle{unsrt}
{\scriptsize

\newpage
In this Supplemental Material we provide the following additional plots:
\begin{enumerate}

    \item
    Different norms and their time averages for the temporally frozen states in $d=2$.
    \item
    In $d=3$, we plot versus $t$, $Q_{0,m}$, $Q_{1,m}$, and $D_{m}$ and their time averages versus $m$, for a representative run.  Weighted averages of the vorticity are defined as
   \begin{equation}
   D_{m}(t)=\|\bom\|^{\alpha_{n,m,d=3}}_{2m}\, \nonumber.     
   \end{equation}
    \end{enumerate}

    \section{Supplemental results for two dimensions}


	\begin{figure*}[!h]
		\begin{tikzpicture}
		\node[anchor=south west,inner sep=0] at (0,0)
		{\includegraphics[width=0.5\linewidth]{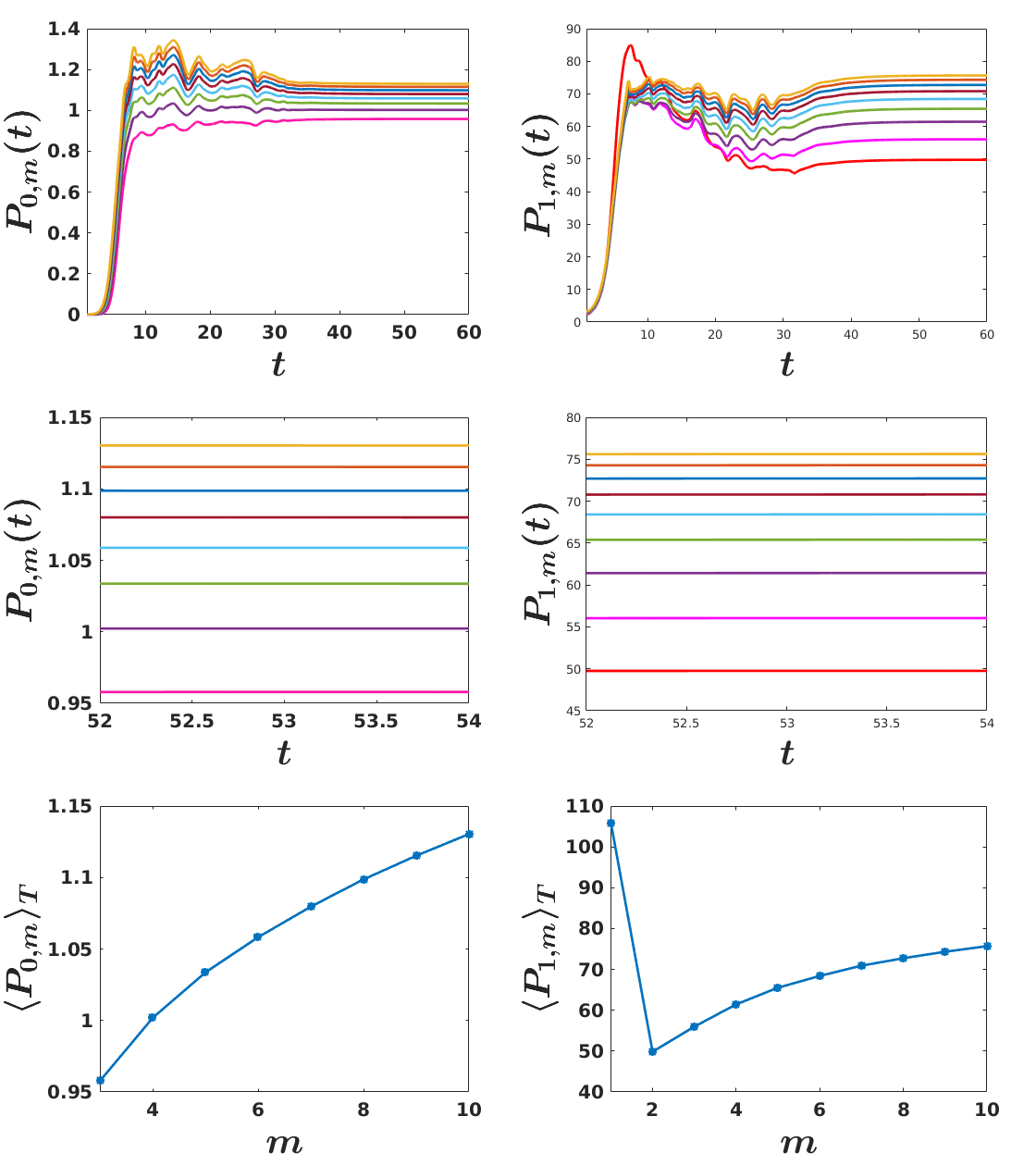}
		\put(-100,225){\large A}
		\put(90,225){\large B}
	
		\includegraphics[width=0.5\linewidth]{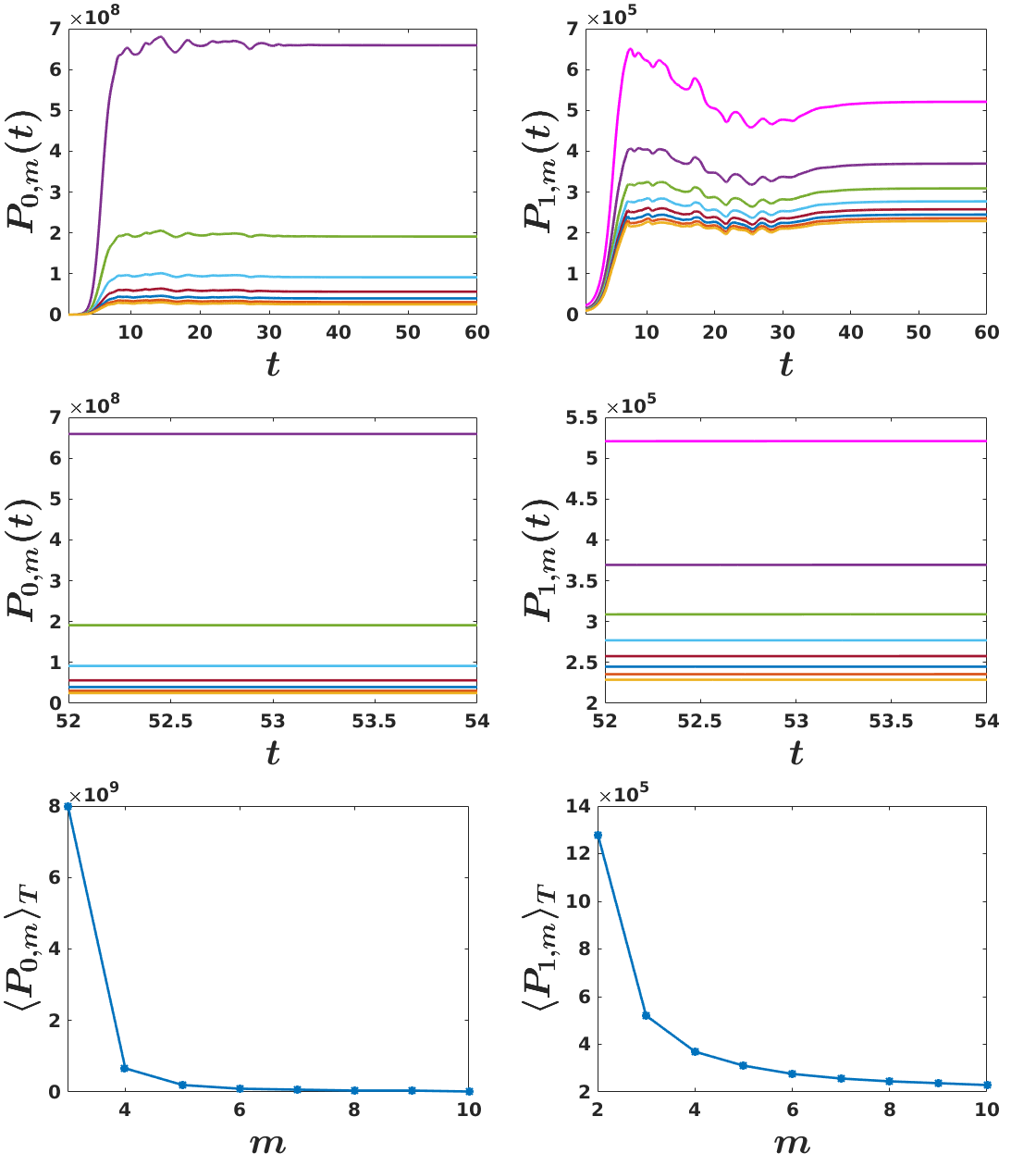}};
		\draw[line width=0.2mm,black, -] (0,0.0) -- (6.7,0.0);
		\draw[line width=0.2mm,black, -] (0,0.0) -- (0.0,7.8);
		\draw[line width=0.2mm,black, -] (6.7,0.0) -- (6.7,7.8);
		\draw[line width=0.2mm,black, -] (0,7.8) -- (6.7,7.8);
		\draw[line width=0.1mm,blue, dashed] (1.825,5.55) -- (0.7,4.9);
		\draw[line width=0.1mm,blue, dashed] (2.0,5.55) -- (3.1,4.9);
		\draw[line width=0.1mm,blue, dashed] (5.2,5.55) -- (4.0,4.9);
		\draw[line width=0.1mm,blue, dashed] (5.3,5.55) -- (6.4,4.9);
		\draw[line width=0.2mm,black, -] (6.75,0) -- (13.6,0);
		\draw[line width=0.2mm,black, -] (6.75,0) -- (6.75,7.8);
		\draw[line width=0.2mm,black, -] (6.75,7.8) --(13.6,7.8);
		\draw[line width=0.2mm,black, -] (13.6,0) --(13.6,7.8);
		\draw[line width=0.1mm,blue, dashed] (8.625,5.55) -- (7.5,4.9);
		\draw[line width=0.1mm,blue, dashed] (8.8,5.55) -- (9.9,4.9);
		\draw[line width=0.1mm,blue, dashed] (11.9,5.55) -- (10.8,4.9);
		\draw[line width=0.1mm,blue, dashed] (12.1,5.55) -- (13.2,4.9);
		\end{tikzpicture}
		
			\caption{(Colour online) Illustrative plots for $U_{0} = \sqrt{\alpha/\beta}$ (panel A) and $U_{0} = \nu/L$ (panel B) for run F7 (see Table 1, main text): First and second rows\,: plots versus $t$ of $P_{0,m}$ and $P_{1,m}$\,; the plots in the second row are expanded versions of small segments of the plots in the first row. Third row\,: Plots versus $m$ of $\left< P_{0,m} \right>_{T}$ and $\left< P_{1,m} \right>_{T}$. Curves for
			$m=2, 3, 4, 5, 6, 7, 8, 9$, and $10$ are drawn in red, pink; violet, green, cyan, maroon, blue, orange, and yellow, respectively.}
				 \label{fig:SM1}
		\end{figure*}
\newpage
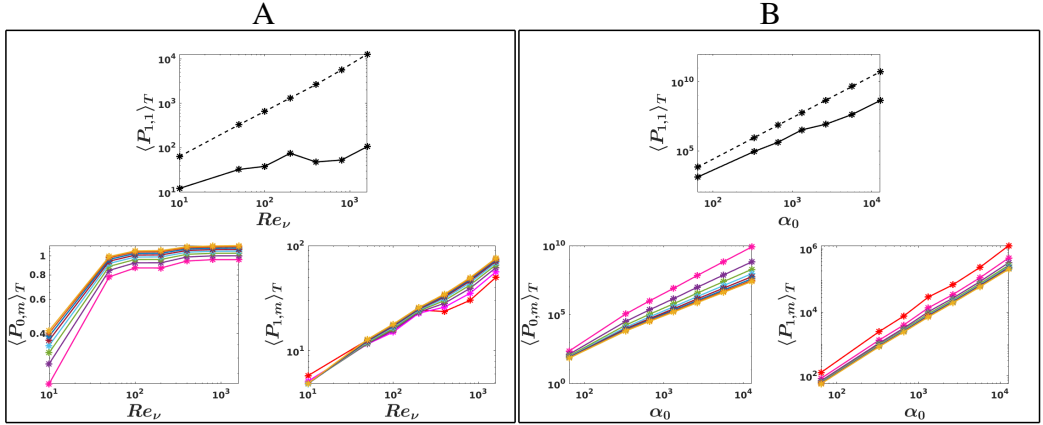
\begin{figure*}[!]
		\begin{tikzpicture}
		\node[anchor=south west,inner sep=0] at (0,0)
		{\includegraphics[width=0.5\linewidth]{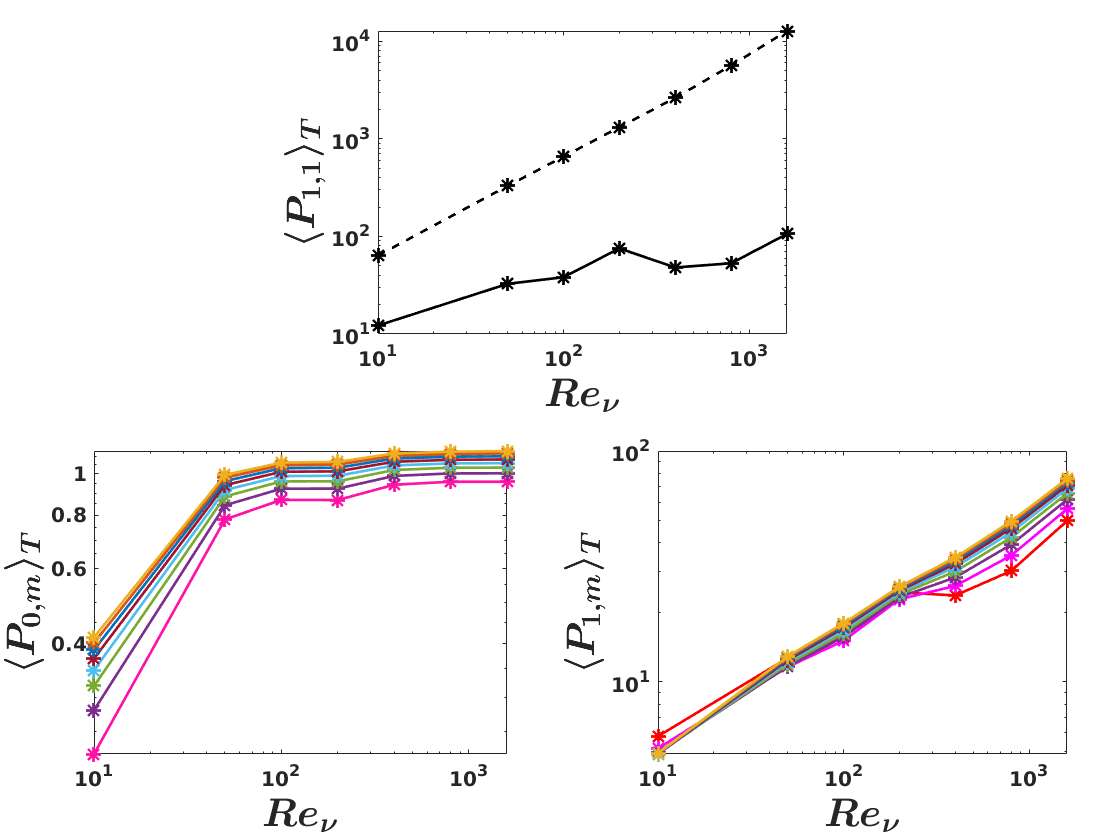}
		\put(-100,150){\large A}
		\put(90,150){\large B}
	
		\includegraphics[width=0.5\linewidth]{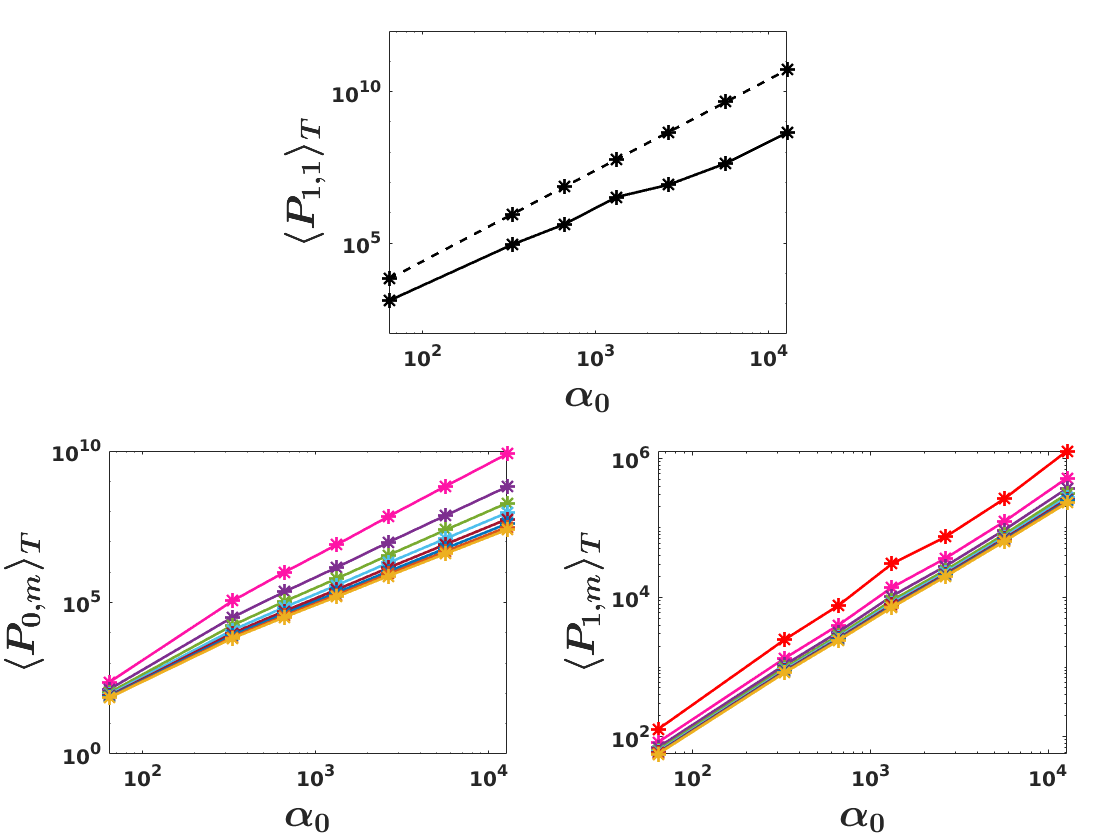}};
		\draw[line width=0.2mm,black, -] (0,0.0) -- (6.7,0.0);
		\draw[line width=0.2mm,black, -] (0,0.0) -- (0.0,5.2);
		\draw[line width=0.2mm,black, -] (6.7,0.0) -- (6.7,5.2);
		\draw[line width=0.2mm,black, -] (0,5.2) -- (6.7,5.2);
	
		\draw[line width=0.2mm,black, -] (6.75,0) -- (13.6,0);
		\draw[line width=0.2mm,black, -] (6.75,0) -- (6.75,5.2);
		\draw[line width=0.2mm,black, -] (6.75,5.2) --(13.6,5.2);
		\draw[line width=0.2mm,black, -] (13.6,0) --(13.6,5.2);
		
		\end{tikzpicture}
		\caption{(Colour online) Illustrative plots for $U_{0} = \sqrt{\alpha/\beta}$ (panel A) and $U_{0} = \nu/L$ (panel B) for $d=2$, runs F1-F7 (see Table 1, main text)\,: First  row\,: plots versus $Re_{\nu}$ (panel A) and $\alpha_{0}$ (panel B) of $\left< P_{1,1} \right> _T$ ( solid black line) and $Re_{\nu}\,\alpha_{0}\,\mathcal{A}_{0}$ (dashed black line). Second row\,: Plots versus $Re_{\nu}$ (panel A) and $\alpha_{0}$ (panel B) of $\left<  P_{0,m}  \right> _T$ and $\left<P_{1,m}\right> _T$. Curves for
			$m=2, 3, 4, 5, 6, 7, 8, 9$, and $10$ are drawn in red, pink; violet, green, cyan, maroon, blue, orange, and yellow, respectively.}\label{fig:SM2}
		\end{figure*}

\medskip
\medskip
\medskip
\section{Supplemental results for three dimensions}

		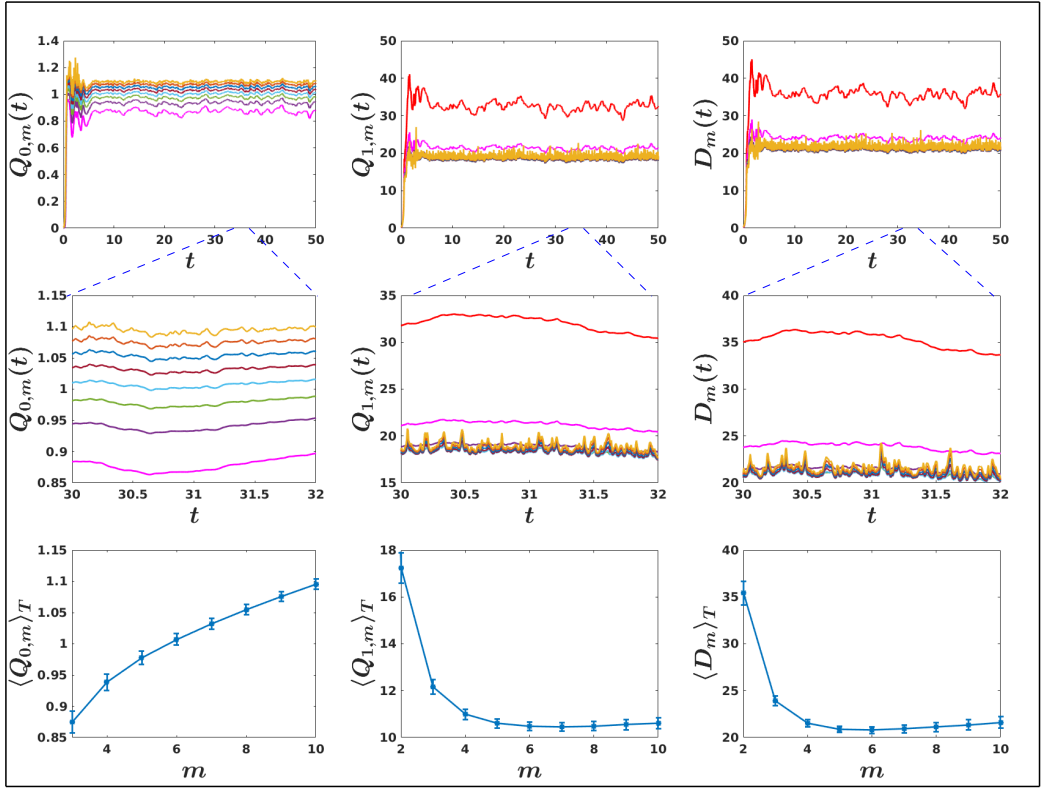
\begin{figure*}[!h]

		\begin{tikzpicture}
		\node[anchor=south west,inner sep=0] at (0,0)
		{
		\includegraphics[width=\linewidth]{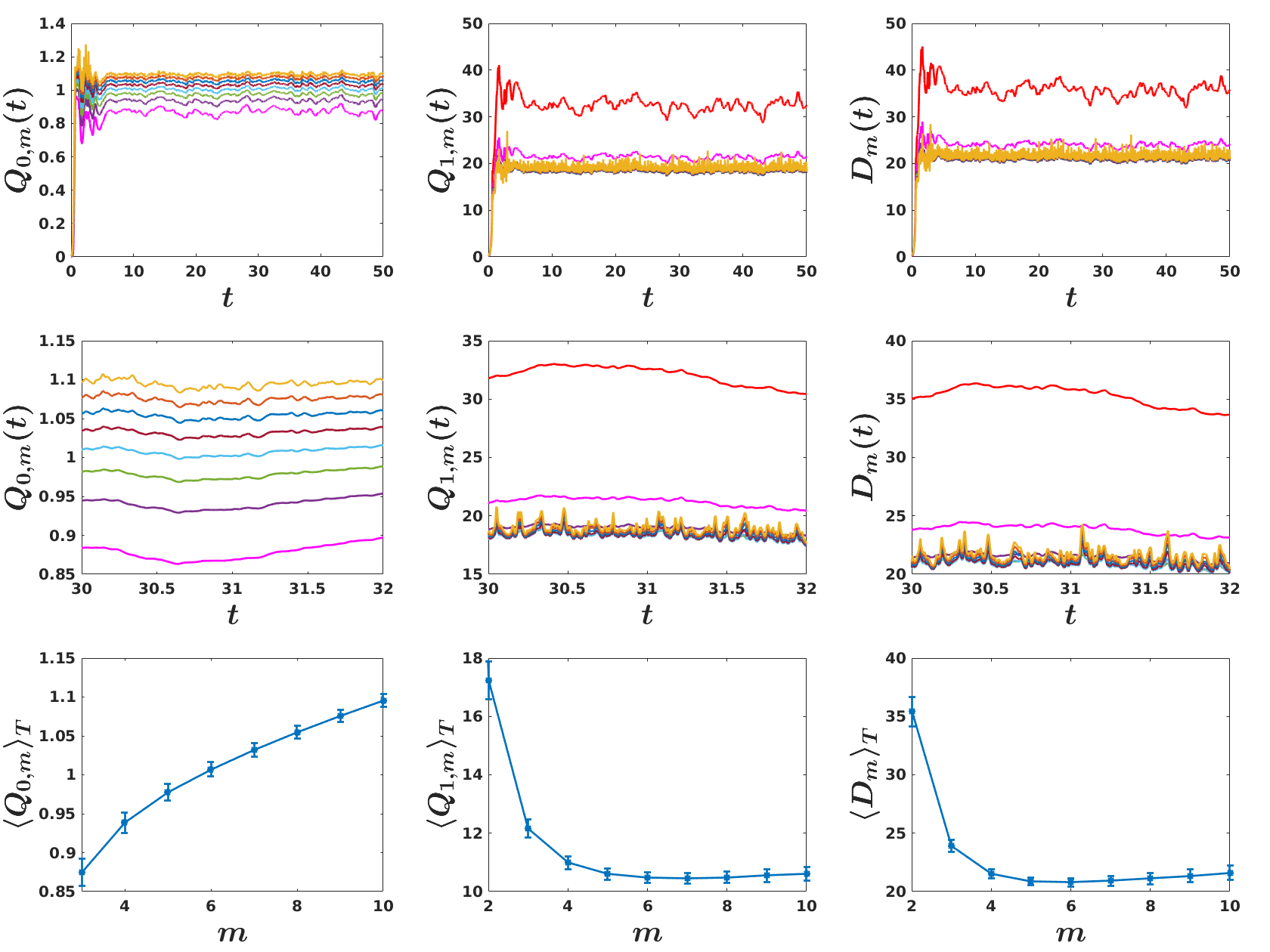}
		\put(-200,310){\large A}
		};
	\draw[line width=0.2mm,black, -] (0,0.0) -- (13.6,0.0);
	\draw[line width=0.2mm,black, -] (0,0.0) -- (0.0,10.4);
	\draw[line width=0.2mm,black, -] (13.6,0.0) -- (13.6,10.4);
	\draw[line width=0.2mm,black, -] (0,10.4) -- (13.6,10.4);
		\draw[line width=0.1mm,blue, dashed] (3.0,7.4) -- (0.8,6.5);
		\draw[line width=0.1mm,blue, dashed] (3.2,7.4) -- (4.1,6.5);
		\draw[line width=0.1mm,blue, dashed] (7.4,7.4) -- (5.3,6.5);
		\draw[line width=0.1mm,blue, dashed] (7.6,7.4) -- (8.5,6.5);
		\draw[line width=0.1mm,blue, dashed] (11.8,7.4) -- (9.7,6.5);
		\draw[line width=0.1mm,blue, dashed] (12.0,7.4) -- (13.0,6.5);
		
		\end{tikzpicture}
		
			\caption{(Colour online) Illustrative plots for $U_{0} = \sqrt{\alpha/\beta}$ (panel A) run B3 (see Table 1, main text): First and second rows\,: plots versus $t$ of $Q_{0,m}$, $Q_{1,m}$ and $D_{m}$\,; the plots in the second row are expanded versions of small segments of the plots in the first row. Third row\,: Plots versus $m$ of $\left< Q_{0,m} \right>_{T}$ and $\left< Q_{1,m} \right>_{T}$ and $\left< D_{m} \right>_{T}$. Curves for
			$m=2, 3, 4, 5, 6, 7, 8, 9$, and $10$ are drawn in red, pink; violet, green, cyan, maroon, blue, orange, and yellow, respectively.}
				\label{fig:SM4}
		\end{figure*}
		\begin{figure*}[!h]

		\begin{tikzpicture}
		\node[anchor=south west,inner sep=0] at (0,0)
		{
	
		\includegraphics[width=\linewidth]{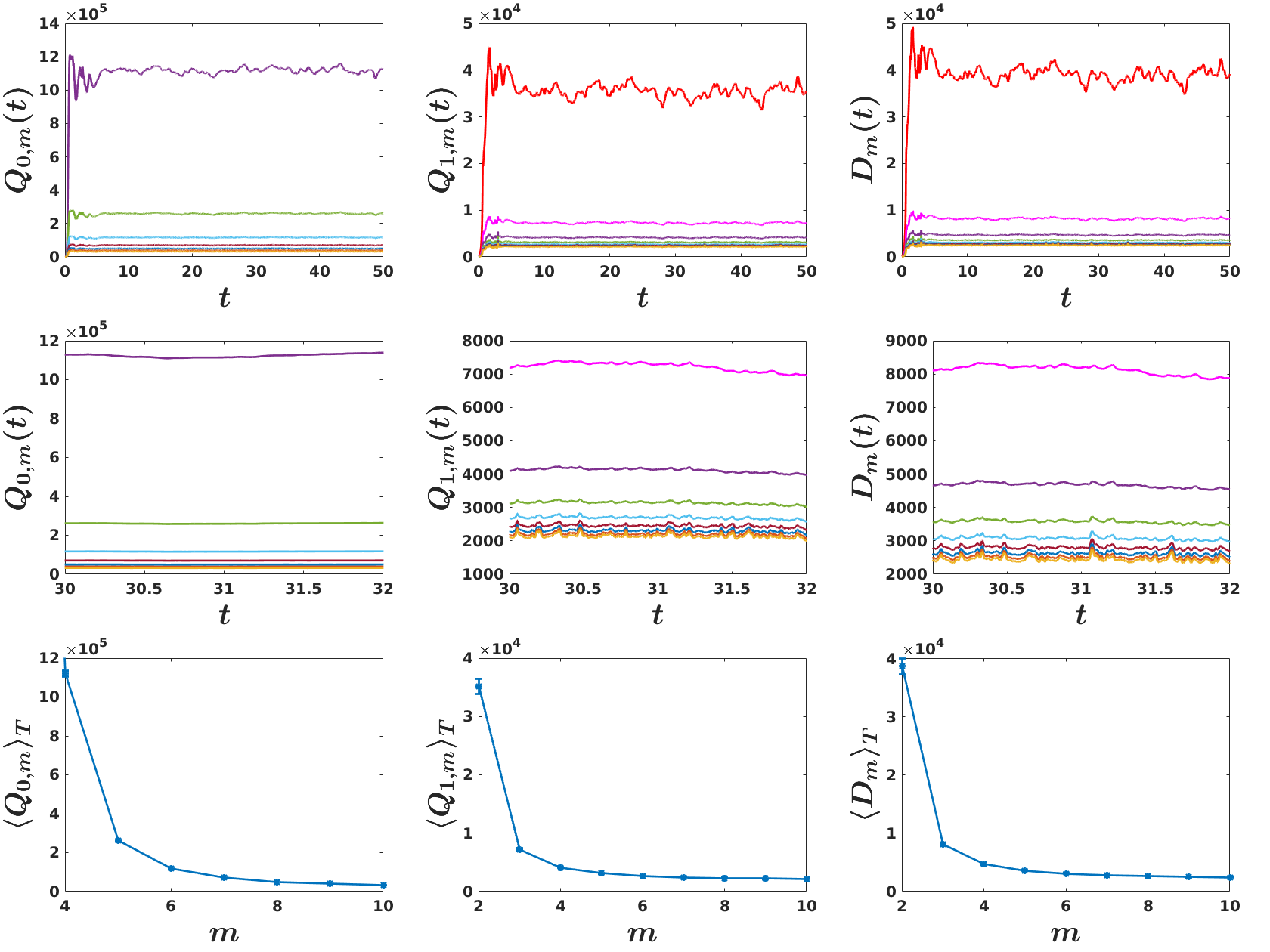}
		\put(-200,310){\large B}
		};
	\draw[line width=0.2mm,black, -] (0,0.0) -- (13.6,0.0);
	\draw[line width=0.2mm,black, -] (0,0.0) -- (0.0,10.4);
	\draw[line width=0.2mm,black, -] (13.6,0.0) -- (13.6,10.4);
	\draw[line width=0.2mm,black, -] (0,10.4) -- (13.6,10.4);
		\draw[line width=0.1mm,blue, dashed] (3.0,7.4) -- (0.8,6.5);
		\draw[line width=0.1mm,blue, dashed] (3.2,7.4) -- (4.1,6.5);
		\draw[line width=0.1mm,blue, dashed] (7.4,7.4) -- (5.3,6.5);
		\draw[line width=0.1mm,blue, dashed] (7.6,7.4) -- (8.5,6.5);
		\draw[line width=0.1mm,blue, dashed] (11.8,7.4) -- (9.7,6.5);
		\draw[line width=0.1mm,blue, dashed] (12.0,7.4) -- (13.0,6.5);
	
		\end{tikzpicture}
		
			\caption{(Colour online) Illustrative plots for $U_{0} = \nu/L$ (panel B) for run B3 (see Table 1, main text): First and second rows\,: plots versus $t$ of $Q_{0,m}$, $Q_{1,m}$ and $D_{m}$\,; the plots in the second row are expanded versions of small segments of the plots in the first row. Third row\,: Plots versus $m$ of $\left< Q_{0,m} \right>_{T}$ and $\left< Q_{1,m} \right>_{T}$ and $\left< D_{m} \right>_{T}$. Curves for
			$m=2, 3, 4, 5, 6, 7, 8, 9$, and $10$ are drawn in red, pink; violet, green, cyan, maroon, blue, orange, and yellow, respectively.}
				\label{fig:SM5}
		\end{figure*}

\end{document}